\documentclass[reqno, a4paper,12pt]{amsart}

\usepackage[capposition=top]{floatrow}
\usepackage{euscript,eufrak,verbatim}
\usepackage{graphicx}
\usepackage[usenames]{color}
\usepackage[colorlinks,linkcolor=red,anchorcolor=blue,citecolor=blue]{hyperref}
\usepackage{amsmath, mathtools}
\usepackage{mathtools}
\usepackage{amsthm}
\usepackage[all]{xy}
\usepackage{amssymb} 
\usepackage{bm}

\usepackage[all]{xy}

\usepackage{mathrsfs}
\usepackage{amscd}

\usepackage{enumitem}

\makeatletter
\numberwithin{equation}{section}

\theoremstyle{plain}
\newtheorem{main theorem}{Main Theorem}
\newtheorem{theorem}{Theorem}[section]
\newtheorem{lemma}[theorem]{Lemma}

\newtheorem{corollary}[theorem]{Corollary}
\newtheorem{proposition}[theorem]{Proposition}
\newtheorem{claim}[theorem]{Claim}

\theoremstyle{definition}
\newtheorem{definition}[theorem]{Definition}
\newtheorem{remark}[theorem]{Remark}
\newtheorem{example}[theorem]{Example}

\newcommand{\spa}{\hspace{1pt}}
\newcommand{\vep}{\varepsilon}
\newcommand{\flo}{\mathscr}
\newcommand{\diam}{\mathrm{diam}}
\newcommand{\tpi}{\widetilde{\pi}}
\DeclarePairedDelimiter{\abs}{\lvert}{\rvert}
\DeclarePairedDelimiter{\ceil}{\lceil}{\rceil}
\newcommand{\norm}[1]{\left\lVert#1\right\rVert}

\begin{document}
\title[Weighted topological pressure revisited]{Weighted topological pressure revisited}

\author[Nima Alibabaei]{Nima Alibabaei}



\maketitle

\begin{abstract}
Feng--Huang (2016) introduced weighted topological entropy and pressure for factor maps between dynamical systems and established its variational principle. Tsukamoto (2022) redefined those invariants quite differently for the simplest case and showed via the variational principle that the two definitions coincide. We generalize Tsukamoto's approach, redefine the weighted topological entropy and pressure for higher dimensions, and prove the variational principle. Our result allows for an elementary calculation of the Hausdorff dimension of affine-invariant sets such as self-affine sponges and certain sofic sets that reside in Euclidean space of arbitrary dimension.
\end{abstract}

\section{Introduction} \label{section: introduction}

\subsection{Dynamical systems and entropy}

Topological pressure and its variational principle have been significant in several fields, including the dimension theory of dynamical systems. Recently, Feng--Huang devised an innovative invariant called weighted topological pressure for factor maps between dynamical systems and proved its variational principle \cite{Feng--Huang}. Their work inspired Tsukamoto to suggest a new definition for this invariant \cite{Tsukamoto}. He also established a variational principle, revealing the nontrivial coincidence of the two definitions. Tsukamoto focused on the simplest case with two dynamical systems.

In this paper, we extend Tsukamoto's definition to the case of an arbitrary number of dynamical systems and prove its variational principle. With our result, we can plainly calculate the Hausdorff dimension of self-affine sponges, a topic studied by Kenyon--Peres \cite{Kenyon--Peres}. Furthermore, we will show in section \ref{section: example: sofic set} that we can determine the Hausdorff dimension of certain sofic sets embedded in higher-dimensional Euclidean space.

We review the basic notions of dynamical systems in this subsection. Refer to the book of Walters \cite{Walters} for the details.

A pair $(X, T)$ is called a \textbf{dynamical system} if $X$ is a compact metrizable space and $T: X \rightarrow X$ is a continuous map. A map $\pi: X \rightarrow Y$ between dynamical systems $(X, T)$ and $(Y, S)$ is said to be a \textbf{factor map} if $\pi$ is a continuous surjection and $\pi \circ T = S \circ \pi$. We sometimes write as $\pi: (X, T) \rightarrow (Y, S)$ to clarify the dynamical systems in question.

For a dynamical system $(X, T)$, denote its \textbf{topological entropy} by $h_{\mathrm{top}}(T)$. Let $P(f)$ be the \textbf{topological pressure} for a continuous function $f: X \rightarrow \mathbb{R}$ (see section \ref{section: weighted topological pressure} for the definition of these quantities). Let $\flo{M}^T(X)$ be the set of $T$-invariant probability measures on $X$ and $h_\mu(T)$ the \textbf{measure-theoretic entropy} for $\mu \in \flo{M}^T(X)$ (see subsection \ref{subsection: measure theoretic entropy}). The variational principle then states that \cite{Dinaburg, Goodman, Goodwyn, Ruelle, Walters75}
\[  P(f) = \sup_{\mu \in \flo{M}^T(X)} \left( h_\mu(T) + \int_X fd\mu \right). \]

\subsection{Background} \label{subsection: background}

We first look at {\it{self-affine sponges}} to understand the background of weighted topological entropy introduced by Feng--Huang. Let $m_1, m_2, \ldots, m_r$ be natural numbers with $m_1 \leq m_2 \leq \cdots \leq m_r$. Consider an endomorphism $T$ on $\mathbb{T}^r = \mathbb{R}^r/\mathbb{Z}^r$ represented by the diagonal matrix $A = \mathrm{diag}(m_1, m_2, \ldots, m_r)$.
For $D \subset \prod_{i=1}^r \{0, 1, \ldots, m_i-1\}$, define
\[ K(T, D) = \left\{ \sum_{n=0}^{\infty} A^{-n}e_n \in \mathbb{T}^r \spa \middle| \spa e_n \in D \right\}. \]
This set is compact and $T$-invariant, i.e., $TK(T, D) = T$.

These sets for $r = 2$ are known as {\it{Bedford-McMullen carpets}} or {\it{self-affine carpets}}. The following figure exhibits a famous example, the case of $D = \{(0,0), (1,1), (0,2)\} \subset \{0, 1\} \times \{0, 1, 2\}$.
\begin{figure}[h!]
\advance\rightskip-2cm
\includegraphics[width=14cm]{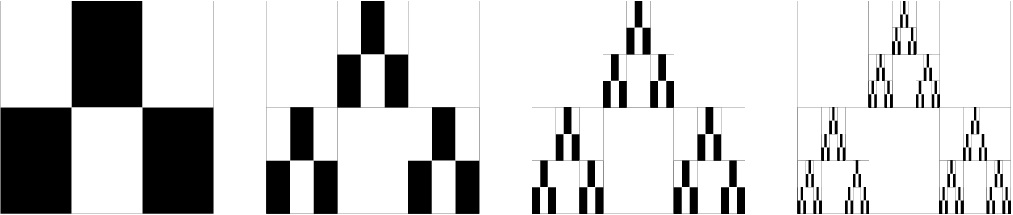}
\vspace{20pt}%
\caption{First four generations of Bedford-McMullen carpet}
\end{figure}
The analysis of these sets is complicated compared to ``self-similar'' sets.
\linebreak
Bedford \cite{Bedford} and McMullen \cite{McMullen} independently studied these sets and showed that, in general, their Hausdorff dimension is strictly smaller than their Minkowski dimension (a.k.a. \hspace{-20pt} Box-counting dimension). The figure above has Hausdorff dimension $\log_2{(1+2^{\log_3{2}})} = 1.349 \cdots$ and Minkowski dimension $1 + \log_3{\frac{3}{2}} = 1.369 \cdots$.

The sets $K(T, D)$ for $r \geq 3$ are called {\it{self-affine sponges}}. Kenyon--Peres \cite{Kenyon--Peres} calculated their Hausdorff dimension for the general case (see Theorem \ref{theorem: self-affine sponges dimension} in this section). In addition, they showed the following variational principle for the Hausdorff dimension of $K(T, D)$;
\begin{equation} \label{equation: Ledrappier-Young formula}
\mathrm{dim}_H K(T, D) = \sup_{\mu \in \flo{M}^T(\mathbb{T}^r)}{ \left\{ 
\frac{1}{\log{m_r}}h_{\mu}(T) + \sum_{i=2}^r \left( \frac{1}{\log{m_{r-i+1}}} - \frac{1}{\log{m_{r-i+2}}} \right) h_{\mu_i}(T_i) \right\}}.
\end{equation}
Here, the endomorphism $T_i$ on $\mathbb{T}^{r-i+1}$ is defined from $A_i = \mathrm{diag}(m_1, m_2, \ldots, m_{r-i+1})$, and $\mu_i$ is defined as the push-forward measure of $\mu$ on $\mathbb{T}^{r-i+1}$ by the projection onto the first $r-i+1$ coordinates. Feng--Huang's definition of weighted topological entropy of $K(T, D)$ equals $\mathrm{dim}_H K(T, D)$ with a proper setting.

\subsection{The original definition of the weighted topological pressure} \label{subsection: the original definition of the weighted topological pressure}

Motivated by the geometry of self-affine sponges described in the previous subsection, Feng--Huang introduced a generalized notion of pressure. Consider dynamical systems $(X_i, \spa T_i)$ ($i=1, \spa 2, \spa \ldots, \spa r$) and factor maps $\pi_i: X_i \rightarrow X_{i+1} \hspace{5pt} (i=1, \spa 2, \spa \ldots \spa , \spa r-1)$:
\begin{equation*}
\xymatrix{
(X_1, T_1) \ar[r]^-{\pi_1} & (X_2, T_2) \ar[r]^-{\pi_2} & \cdots \ar[r]^-{\pi_{r-1}} & (X_r, T_r)}.
\end{equation*}
We refer to this as a {\textbf{sequence of dynamical systems}}. Let the weight $\boldsymbol{w} = (w_1, w_2, \ldots, w_r)$ with $w_1 > 0$ and $w_i \geq 0$ for $i \geq 2$. Feng--Huang \cite{Feng--Huang} ingeniously defined the $\boldsymbol{w}$-weighted topological pressure $P^{\boldsymbol{w}}_\mathrm{FH}(f)$ for a continuous function $f:X_1 \rightarrow \mathbb{R}$ and established the variational principle \cite[Theorem 1.4]{Feng--Huang}
\begin{equation} \label{equation: FH variational principle}
P^{\boldsymbol{w}}_{\mathrm{FH}}(f) = \sup_{\mu \in \mathscr{M}^{T_1}(X_1)} \left( \sum_{i=1}^r w_i h_{{\pi^{(i-1)}}_*\mu} (T_i) + w_1 \int_{X_1} f d\mu \right).
\end{equation}
Here $\pi^{(i)}$ is defined by
\begin{gather*}
\pi^{(0)} = \mathrm{id}_{X_1}: X_1 \to X_1, \\
\pi^{(i)} = \pi_i \circ \pi_{i-1} \circ \cdots \circ \pi_1: X_1 \to X_{i+1},
\end{gather*}
and ${\pi^{(i-1)}}_* \mu$ is the push-forward measure of $\mu$ by $\pi^{(i-1)}$ on $X_i$. The $\boldsymbol{w}$-weighted topological entropy $h^{\boldsymbol{w}}_{\mathrm{top}}(T_1)$ is the value of $P^{\boldsymbol{w}}_\mathrm{FH}(f)$ when $f \equiv 0$. In this case, \eqref{equation: FH variational principle} becomes
\begin{equation} \label{equation: FH variational principle for entropy}
h^{\boldsymbol{w}}_{\mathrm{top}}(T_1) = \sup_{\mu \in \flo{M}^{T_1}(X_1)} \left( \sum_{i=1}^r w_i h_{{\pi^{(i-1)}}_*\mu} (T_i) \right).
\end{equation}

We will explain here Feng--Huang's method of defining $h^{\boldsymbol{w}}_{\mathrm{top}}(T_1)$. For the definition of $P^{\boldsymbol{w}}_\mathrm{FH}(f)$, see their original paper \cite{Feng--Huang}.

Let $n$ be a natural number and $\vep$ a positive number. Let $d^{(i)}$ be a metric on $X_i$. For $x \in X_1$, define the \textbf{$\boldsymbol{n}$-th $\boldsymbol{w}$-weighted Bowen ball of radius $\boldsymbol{\vep}$ centered at $\boldsymbol{x}$} by
\begin{equation*}
B^{\boldsymbol{w}}_n(x, \vep) = \left\{ y \in X_1 \spa \middle|
\begin{array}{l}
\text{$d^{(i)} \! \left( T^j_i(\pi^{(i-1)}(x)), T^j_i(\pi^{(i-1)}(y)) \right) < \vep$ for every} \\[2pt]
\text{$0 \leq j \leq \ceil{(w_1 + \cdots + w_i)n}$ and $1 \leq i \leq k$.}
\end{array}
\right\}.
\end{equation*} 

Consider $\Gamma = \{ B^{\boldsymbol{w}}_{n_j}(x_j, \vep) \}_j$, an at-most countable cover of $X_1$ by Bowen balls. Let $n(\Gamma) = \min_j n_j$. For $s \geq 0$ and $N \in \mathbb{N}$, let
\[ \Lambda^{\boldsymbol{w}, s}_{N, \vep} = \inf \left\{ \sum_j e^{-sn_j} \spa \middle| \spa \text{ $\Gamma = \{ B^{\boldsymbol{w}}_{n_j}(x_j, \vep) \}_j$ covers $X_1$ and $n(\Gamma) \geq N$} \right\}. \]
This quantity is non-decreasing as $N \to \infty$. The following limit hence exists:
\[ \Lambda^{\boldsymbol{w}, s}_{\vep} = \lim_{N \to \infty} \Lambda^{\boldsymbol{w}, s}_{N, \vep}. \]
There is a value of $s$ where $\Lambda^{\boldsymbol{w}, s}_{\vep}$ jumps from $\infty$ to $0$, which we will denote by $h^{\boldsymbol{w}}_{\mathrm{top}}(T_1, \vep)$:
\begin{equation*}
\Lambda^{\boldsymbol{w}, s}_{\vep} = \left\{
\begin{array}{ll}
\infty & (s < h^{\boldsymbol{w}}_{\mathrm{top}}(T_1, \vep)) \\
0 & (s > h^{\boldsymbol{w}}_{\mathrm{top}}(T_1, \vep))
\end{array}
\right..
\end{equation*}
The value $h^{\boldsymbol{w}}_{\mathrm{top}}(T_1, \vep)$ is non-decreasing as $\vep \to 0$. Therefore, we can define the $\boldsymbol{w}$-weighted topological entropy $h^{\boldsymbol{w}}_{\mathrm{top}}(T_1)$ by
\[ h^{\boldsymbol{w}}_{\mathrm{top}}(T_1) = \lim_{\vep \to 0} h^{\boldsymbol{w}}_{\mathrm{top}}(T_1, \vep). \]

An important point about this definition is that in some dynamical systems, such as self-affine sponges, the quantity $h^{\boldsymbol{w}}_{\mathrm{top}}(T_1)$ is directly related to the Hausdorff dimension of $X_1$.

\begin{example} \label{example: self affine sponges 1}
Consider the self-affine sponges introduced in subsection \ref{subsection: background}. Define $p_i: \mathbb{T}^{r-i+1} \rightarrow \mathbb{T}^{r-i}$ by
\[ p_i(x_1, x_2, \ldots, x_{r-i}, x_{r-i+1}) = (x_1, x_2, \ldots, x_{r-i}). \]
Let $X_1 = K(T, D)$, $X_i = p_{i-1} \circ p_i \circ \cdots \circ p_1(X_1)$, and $T_i: X_i \rightarrow X_i$ be the endomorphism defined from $A_i = \mathrm{diag}(m_1, m_2, \ldots, m_{r-i+1})$. Define the factor maps $\pi_i: X_i \rightarrow X_{i+1}$ as the restrictions of $p_i$. Let
\begin{equation} \label{equation: weight for self-affine sponges}
\boldsymbol{w} = \left( \frac{\log{m_1}}{\log{m_r}}, \quad \frac{\log{m_1}}{\log{m_{r-1}}} - \frac{\log{m_1}}{\log{m_r}}, \spa \ldots \spa , \quad \frac{\log{m_1}}{\log{m_2}} - \frac{\log{m_1}}{\log{m_3}}, \quad 1 - \frac{\log{m_1}}{\log{m_2}} \right).
\end{equation}
Then $n$-th $\boldsymbol{w}$-weighted Bowen ball is approximately a square with a side length of $\vep m_1^{-n}$. Therefore,
\begin{equation} \label{equation: self affine sponges through FH}
\mathrm{dim}_H K(T, D) = \frac{h^{\boldsymbol{w}}_{\mathrm{top}}(T_1)}{\log{m_1}}.
\end{equation}
\end{example}

\subsection{Tsukamoto's approach and its extension} \label{subsection: Tsukamoto's approach and its extension}
Following the work of Feng--Huang \cite{Feng--Huang} described in the previous subsection, Tsukamoto \cite{Tsukamoto} published an intriguing approach to these invariants. There, he gave a new definition of the weighted topological pressure for two dynamical systems and a factor map:
\[ \xymatrix{
(X_1, T_1) \ar[r]^-{\pi} & (X_2, T_2).} \]
He then proved the variational principle using his definition, showing the surprising coincidence of the two definitions. His expression of weighted topological entropy allowed for relatively easy calculations for sets like self-affine carpets.

We will extend Tsukamoto's idea, redefine the weighted topological pressure for an arbitrary length of a sequence of dynamical systems, and establish the variational principle. Here we will explain our definition in the case $f \equiv 0$. See section \ref{section: weighted topological pressure} for the general setting. We will not introduce Tsukamoto's definition since it is obtained by letting $r=2$ in the following argument.

Let ${\boldsymbol{a}}=(a_1, \spa a_2, \spa \cdots, a_{r-1})$ with $0 \leq a_i \leq 1$ for each $i$. Let $N$ be a natural number and $\vep$ a positive number. We define a new metric $d^{(i)}_N$ on $X_i$ by
\[ d^{(i)}_N(x_1, \spa x_2) = \max_{0\leq n < N} d^{(i)}({T_i}^n x_1, \spa {T_i}^n x_2). \]
For $\Omega \subset X_1$, we define
\begin{align*}
\#^{\boldsymbol{a}}_1(\Omega, \spa N, \spa \vep)
&= \min \left\{ \spa n \in \mathbb{N} \spa \middle|
\begin{array}{l}
\text{There exists an open cover $\{U_j\}_{j=1}^n$ of $\Omega$} \\ 
\text{with $\diam(U_j, \spa d_N^{(1)}) < \vep$ for all $ 1 \spa \leq j \spa \leq n$}
\end{array}
\right\}.
\end{align*} 
Let $\Omega \subset X_{i+1}$. If $\#^{\boldsymbol{a}}_i$ is already defined, let
\begin{flalign*}
& \#^{\boldsymbol{a}}_{i+1}(\Omega, \spa N, \spa \vep) &
\end{flalign*} \\[-35pt]
\begin{align*}
&= \min \left\{ \spa \sum_{j=1}^n \Big( \#^{\boldsymbol{a}}_i(\pi_i^{-1}(U_j), \spa N, \spa \vep) \Big)^{a_i} \spa \middle|
\begin{array}{l}
\text{$n \in \mathbb{N}$, $\{U_j\}_{j=1}^n$ is an open cover of $\Omega$} \\ 
\text{with $\diam(U_j, \spa d_N^{(i+1)}) < \vep$ for all $ 1 \spa \leq j \spa \leq n$}
\end{array}
\right\}.
\end{align*} 
We define the \textbf{topological entropy of ${\boldsymbol{a}}$-exponent} $h^{\boldsymbol{a}}(\boldsymbol{T})$, where $\boldsymbol{T} = (T_i)_i$, by \\[1pt]
\begin{equation*}
h^{\boldsymbol{a}}(\boldsymbol{T}) = \lim_{\vep \to 0} \left( \lim_{N \to \infty} \frac{\log{\#^{\boldsymbol{a}}_r(X_r, \spa N, \spa \vep)}}{N} \right).
\end{equation*} \\[1pt]
This limit exists since $\log{\#^{\boldsymbol{a}}_r(X_r, \spa N, \spa \vep)}$ is sub-additive in $N$ and non-decreasing as $\vep$ tends to $0$.

From $\boldsymbol{a}$, define $\boldsymbol{w_a} = (w_1, \spa\cdots, \spa w_r)$ by
\begin{eqnarray*}
\left\{
\begin{array}{l}
w_1 = a_1 a_2 a_3 \cdots a_{r-1}\\
w_2 = (1-a_1) a_2 a_3 \cdots a_{r-1} \\
w_3 = (1-a_2) a_3 \cdots a_{r-1} \\
\hspace{50pt} \vdots \\
w_{r-1} = (1-a_{r-2}) a_{r-1} \\
w_r = 1- a_{r-1}
\end{array}
\right..
\end{eqnarray*}
Then our main result Theorem \ref{theorem: main theorem} below yields

\begin{theorem} \label{theorem: main theorem for entropy}
For ${\boldsymbol{a}}=(a_1, \spa a_2, \spa \cdots, a_{r-1})$ with $0 \leq a_i \leq 1$ for each $i$,
\begin{equation} \label{equation: main result for f=0}
h^{\boldsymbol{a}}(\boldsymbol{T}) = \sup_{\mu \in \flo{M}^{T_1}(X_1)} \left( \sum_{i=1}^r w_i h_{{\pi^{(i-1)}}_*\mu} (T_i) \right).
\end{equation}
\end{theorem}
The strategy of the proof is adopted from Tsukamoto's paper. However, there are some additional difficulties. Let $h^{\boldsymbol{a}}_{\mathrm{var}}(\boldsymbol{T})$ be the right-hand side of \eqref{equation: main result for f=0}. We use the ``zero-dimensional trick'' for proving $h^{\boldsymbol{a}}(\boldsymbol{T}) \leq h^{\boldsymbol{a}}_{\mathrm{var}}(\boldsymbol{T})$, meaning we reduce the proof to the case where all dynamical systems are zero-dimensional. Merely taking a zero-dimensional extension for each $X_i$ does not work. Therefore we realize this by taking step by step an extension of the whole sequence of dynamical systems (see subsection \ref{subsection: zero-dimensional principal extension}). Then we show $h^{\boldsymbol{a}}(\boldsymbol{T}) \leq h^{\boldsymbol{a}}_{\mathrm{var}}(\boldsymbol{T})$ by using an appropriate measure, the definition of which is quite sophisticated (see $\sigma_N$ in the proof of Theorem \ref{theorem: first half of the main theorem}).
In proving $h^{\boldsymbol{a}}(\boldsymbol{T}) \geq h^{\boldsymbol{a}}_{\mathrm{var}}(\boldsymbol{T})$, the zero-dimensional trick can not be utilized. The proof, therefore, requires a detailed estimation of these values for arbitrary covers, which is more complicated than the original argument in \cite{Tsukamoto}.

Theorem \ref{theorem: main theorem for entropy} and Feng--Huang's version of variational principle \eqref{equation: FH variational principle for entropy} yield
\begin{corollary} \label{corollary: deep result}
For ${\boldsymbol{a}}=(a_1, \spa a_2, \spa \cdots, a_{r-1})$ with $0 < a_i \leq 1$ for each $i$,
\begin{equation*}
h^{\boldsymbol{a}}(\boldsymbol{T}) = h^{\boldsymbol{w_a}}_{\mathrm{top}}(T_1).
\end{equation*}
\end{corollary}
This corollary is rather profound, connecting the two seemingly different quantities. We can calculate the Hausdorff dimension of certain self-affine sets using this result, as seen in the following example and section \ref{section: example: sofic set}.

\begin{example} \label{example: self-affine sponges two}
Let us take another look at self-affine sponges. Kenyon--Peres \cite[Theorem 1.2]{Kenyon--Peres} calculated their Hausdorff dimension as follows (recall that $m_1 \leq m_2 \leq \cdots \leq m_r$).

\begin{theorem} \label{theorem: self-affine sponges dimension}
Define a sequence of real numbers $(Z_j)_j$ as follows. Let $Z_r$ be the indicator of $D$, namely, $Z_r(i_1, \ldots, i_r) = 1$ if $(i_1, \ldots, i_r) \in D$ and $0$ otherwise. Define $Z_{r-1}$ by
\[ Z_{r-1}(i_1, \ldots, i_{r-1}) = \sum_{i_r = 0}^{m_r-1} Z_r(i_1, \ldots, i_{r-1}, i_r). \]
More generally, if $Z_{j+1}$ is already defined, let
\[ Z_j(i_1, \ldots, i_j) = \sum_{i_{j+1} = 0}^{m_{j+1}-1} Z_{j+1}(i_1, \ldots, i_j, i_{j+1})^{\log{m_{j+1}}/\log{m_{j+2}}}. \]
Then
\[ \mathrm{dim}_H K(T, D) = \frac{\log{Z_0}}{\log{m_1}}. \]
\end{theorem}

We can prove this result fairly elementary by Corollary \ref{corollary: deep result} without requiring measure theory on the surface.
Set $a_i = \log_{m_{r-i+1}} m_{r-i}$ for each $i$, then $\boldsymbol{w_a}$ equals $\boldsymbol{w}$ in \eqref{equation: weight for self-affine sponges}. Combining \eqref{equation: self affine sponges through FH} and Corollary \ref{corollary: deep result},
\[ \mathrm{dim}_H K(T, D) = \frac{h^{\boldsymbol{w_a}}_{\mathrm{top}}(T_1)}{\log{m_1}} = \frac{h^{\boldsymbol{a}}(\boldsymbol{T})}{\log{m_1}}. \]
Hence, we need to show the following claim.

\begin{claim}
We have
\[ h^{\boldsymbol{a}}(\boldsymbol{T}) = \log{Z_0}. \]
\end{claim}

\begin{proof}
Observe first that taking the infimum over closed covers instead of open ones in the definition of $h^{\boldsymbol{a}}(\boldsymbol{T})$ does not change its value. Define a metric $d^{(i)}$ on each $X_i$ by
\[ d^{(i)} (x, y) = \min_{n \in {\mathbb{Z}}^{r-i+1}} \spa \abs{x-y-n}. \]
Let
\[ D_j = \{ (e_1, \ldots, e_j) \spa | \spa \text{ there are $e_{j+1}, \ldots, e_r$ with $(e_1, \ldots, e_r) \in D$} \}. \]
Define $p_i: D_{r-i+1} \rightarrow D_{r-i}$ by $p_i(e_1, \ldots, e_{r-i+1}) = (e_1, \ldots, e_{r-i})$. Fix $0 < \vep < \frac{1}{m_r}$ and take a natural number $n$ with $m_1^{-n} < \vep$. Fix a natural number $N$ and let $\psi_i: D_{r-i+1}^{N+n} \rightarrow D_{r-i}^{N+n}$ be the product map of $p_i$, i.e., $\psi_i(v_1, \ldots, v_{N+n}) = (p_i(v_1), \ldots, p_i(v_{N+n}))$.

For $x \in D_{r-i+1}^{N+n}$, define (recall that $A_i = \mathrm{diag}(m_1, m_2, \ldots, m_{r-i+1})$)
\[ U^{(i)}_x = \left\{ \sum_{k=0}^{\infty} A_i^{-k} e_k \in X_i \spa \middle| \spa \text{$e_k \in D_{r-i+1}$ for each $k$ and $(e_1, \dots, e_{N+n}) = x$} \right\}. \]
Then $\{U^{(i)}_x\}_{x \in D^{N+n}_{r-i+1}}$ is a closed cover of $X_i$ with $\diam(U^{(i)}_x, d^{(i)}_N) < \vep$. For $x, y \in D^{N+n}_{r-i+1}$, we write $x \backsim y$ if and only if $U^{(i)}_x \cap U^{(i)}_y \ne \varnothing$. We have for any $i$ and $x \in D_{r-i}^{N+n}$
\[ \pi_i^{-1}(U^{(i+1)}_x) \subset \bigcup_{\substack{x' \in D_{r-i}^{N+n} \\ x' \backsim x}} \hspace{3pt} \bigcup_{y \in {\psi_i}^{-1}(x')} U^{(i)}_y.  \]
Notice that for each $x \in D_{r-i}^{N+n}$, the number of $x' \in D_{r-i}^{N+n}$ with $x' \backsim x$ is not more than $3^r$. Therefore, for every $v = (v_1^{(1)}, \ldots, v_{N+n}^{(1)}) \in D_{r-1}^{N+n}$, there are $(v_1^{(k)}, \ldots, v_{N+n}^{(k)}) \in D_{r-1}^{N+n}$,  $k = 2, 3, \ldots, L$, and $L \leq 3^r$, with
\begin{equation*}
\#^{\boldsymbol{a}}_1(\pi_1^{-1}(U^{(2)}_v), \spa N, \spa \vep)
\leq \sum_{k=1}^L Z_{r-1}(v_1^{(k)})\cdots Z_{r-1}(v_{N+n}^{(k)}).
\end{equation*}
We inductively continue while considering that the multiplicity is at most $3^r$ and obtain
\begin{flalign*}
& \#^{\boldsymbol{a}}_r(X_r, \spa N, \spa \vep) &
\end{flalign*}
\begin{align*}
&\leq
\begin{multlined}[t][16cm]
3^{r(r-1)} \sum_{x_1 \in D^{N+n}_1}\left( \sum_{x_2 \in {\psi_2}^{-1}(x_1)} \left( \cdots \left( \sum_{x_{r-2} \in {\psi_{r-2}}^{-1}(x_{r-3})} \rule{0cm}{1.1cm} \right. \right. \right. \\
\left. \left. \left. \left( \sum_{\substack{(v_1, \ldots, v_{N+n}) \in {\psi_{r-1}}^{-1}(x_{r-2}) \\ v_j \in D_{r-1} \text{ for each $j$} }} \big( Z_{r-1}(v_1)\cdots Z_{r-1}(v_{N+n}) \big)^{a_1} \right)^{a_2} \right)^{a_3} \cdots \right)^{a_{r-2}}\right)^{a_{r-1}}
\end{multlined} \\
&= 3^{r(r-1)} \left\{ \sum_{x_1 \in D_1} \left( \sum_{x_2 \in p_2^{-1}(x_1)} \left( \cdots \left( \sum_{x_{r-1} \in p_{r-1}^{-1}(x_{r-2})} {Z_{r-1}(x_1, \ldots, x_{r-1})}^{a_1} \right)^{a_2} \cdots \right)^{a_{r-2}}\right)^{a_{r-1}} \right\}^{N+n} \\[5pt]
&= 3^{r(r-1)} {Z_0}^{N+n}.
\end{align*}

\pagebreak
Therefore,
\[ h^{\boldsymbol{a}}(\boldsymbol{T}) = \lim_{\vep \to 0} \left( \lim_{N \to \infty} \frac{\log{\#^{\boldsymbol{a}}_r(X_r, \spa N, \spa \vep)}}{N} \right)
\leq \log{Z_0}. \]

Next, we prove $h^{\boldsymbol{a}}(\boldsymbol{T}) \geq \log{Z_0}$. We fix $0 < \vep < \frac{1}{m_r}$ and utilize $\vep$-separated sets. Take and fix $\boldsymbol{s} = (t_1, \ldots, t_r) \in D$, and set $\boldsymbol{s}_i = (t_1, \ldots, t_{r-i+1})$. Fix a natural number $N$ and let $\psi_i: D_{r-i+1}^N \rightarrow D_{r-i}^N$ be the product map of $p_i$ as in the previous definition. Define
\[ Q_i = \left\{ \sum_{k=1}^N {A_i}^{-k} e_k + \sum_{k=N+1}^{\infty} {A_i}^{-k} \boldsymbol{s}_i \in X_i \spa \middle| \spa e_1, \ldots, e_N \in D_{r-i+1} \right\}. \]
Then $Q_i$ is an $\vep$-separated set with respect to the metric $d^{(i)}_N$ on $X_i$. Consider an arbitrary open cover $\flo{F}^{(i)}$ of $X_i$ for each $i$ with the following properties (this $(\flo{F}^{(i)})_i$ is defined as \textbf{a chain of open ($N$, $\vep$)-covers} of $(X_i)_i$ in Definition \ref{definition: chain of covers}).
\begin{enumerate}
\item[(1)] For every $i$ and $V \in \flo{F}^{(i)}$, we have $\diam(V, d^{(i)}_N) < \vep$.
\item[(2)] For each $1 \leq i \leq r-1$ and $U \in \flo{F}^{(i+1)}$, there is $\flo{F}^{(i)}(U) \subset \flo{F}^{(i)}$ such that
\[ \pi_i^{-1}(U) \subset \bigcup \flo{F}^{(i)}(U) \]
and
\[ \flo{F}^{(i)} = \bigcup_{U \in \flo{F}^{(i+1)}} \flo{F}^{(i)}(U). \]
\end{enumerate}
We have $\#(V \cap Q_i ) \leq 1$ for each $V \in \flo{F}^{(i)}$ by (1). Let $(e^{(2)}_1, e^{(2)}_2, \cdots, e^{(2)}_N) \in D_{r-1}^N$ and suppose $U \in \flo{F}^{(2)}$ satisfies
\[ \sum_{k=1}^N {A_2}^{-k} e^{(2)}_k + \sum_{k=N+1}^{\infty} {A_2}^{-k} \boldsymbol{s}_2 \in U \cap Q_2. \]
Then $\pi_1^{-1}(U)$ contains at least $Z_{r-1}(e^{(2)}_1)\cdots Z_{r-1}(e^{(2)}_N)$ points of $Q_1$. Hence,
\[ \#^{\boldsymbol{a}}_1(\pi_1^{-1}(U), \spa N, \spa \vep)
\geq Z_{r-1}(e^{(2)}_1)\cdots Z_{r-1}(e^{(2)}_N). \]
We continue this reasoning inductively and get
\pagebreak
\begin{flalign*}
& \#^{\boldsymbol{a}}_r(X_r, \spa N, \spa \vep) &
\end{flalign*}
\begin{align*}
&\geq
\begin{multlined}[t][16cm]
\sum_{e^{(1)} \in D^N_1}\left( \sum_{e^{(2)} \in {\psi_2}^{-1}(e^{(1)})} \left( \cdots \left( \sum_{e^{(r-2)} \in {\psi_{r-2}}^{-1}(e^{(r-3)})} \rule{0cm}{1.1cm} \right. \right. \right. \\
\left. \left. \left. \left( \sum_{\substack{(e^{(2)}_1, \ldots, e^{(2)}_N) \in {\psi_{r-1}}^{-1}(e^{(3)}) \\ e^{(2)}_j \in D_{r-1} \text{ for each $j$} }} \big( Z_{r-1}(e^{(2)}_1)\cdots Z_{r-1}(e^{(2)}_N) \big)^{a_1} \right)^{a_2} \right)^{a_3} \cdots \right)^{a_{r-2}}\right)^{a_{r-1}}
\end{multlined} \\
&= \left\{ \sum_{x_1 \in D_1} \left( \sum_{x_2 \in p_2^{-1}(x_1)} \left( \cdots \left( \sum_{x_{r-1} \in p_{r-1}^{-1}(x_{r-2})} {Z_{r-1}(x_1, \ldots, x_{r-1})}^{a_1} \right)^{a_2} \cdots \right)^{a_{r-2}}\right)^{a_{r-1}} \right\}^N \\[5pt]
&= {Z_0}^N.
\end{align*}
This implies
\[ h^{\boldsymbol{a}}(\boldsymbol{T}) \geq \log{Z_0}. \]

We conclude that
\[ h^{\boldsymbol{a}}(\boldsymbol{T}) = \log{Z_0}. \]

\end{proof}
\end{example}

We would like to mention the work of Barral and Feng \cite{Barral--Feng, Feng}, and of Yayama \cite{Yayama}. These papers independently studied the related invariants when $(X, T)$ and $(Y, S)$ are subshifts over finite alphabets.

\section{Weighted topological pressure} \label{section: weighted topological pressure}

Here, we introduce the generalized, new definition of weighted topological pressure. Let $(X_i, \spa T_i)$ ($i=1, \spa 2, \spa \ldots, \spa r$) be dynamical systems and $\pi_i: X_i \rightarrow X_{i+1} \hspace{5pt} (i=1, \spa 2, \spa \ldots \spa , \spa r-1)$ factor maps. For a continuous function $f: X_1 \to \mathbb{R}$ and a natural number $N$, set
\[
S_N f (x) = f(x) + f(T_1 x) + f(T_1^2 x) + \cdots + f(T_1^{N-1}x).
\]
Let $d^{(i)}$ be a metric on $X_i$. Recall that we defined a new metric $d^{(i)}_N$ on $X_i$ by
\[ d^{(i)}_N(x_1, \spa x_2) = \max_{0\leq n < N} d^{(i)}({T_i}^n x_1, \spa {T_i}^n x_2). \]
We may write these as $S_N^{T_1}f$ or $d^{T_i}_N$ to clarify the maps $T_1$ and $T_i$ in the definitions above.

Let ${\boldsymbol{a}}=(a_1, \spa a_2, \spa \cdots, a_{r-1})$ with $0 \leq a_i \leq 1$ for each $i$ and $\vep$ a positive number. For $\Omega \subset X_1$, we define
\pagebreak
\begin{flalign*}
& P^{\boldsymbol{a}}_1(\Omega, \spa f, \spa N, \spa \vep) &
\end{flalign*} \\[-35pt]
\begin{align*}
&= \inf \left\{ \spa \sum_{j=1}^n \exp \left( \sup_{U_j} S_N f \right) \spa \middle|
\begin{array}{l}
\text{$n \in \mathbb{N}$, $\{U_j\}_{j=1}^n$ is an open cover of $\Omega$} \\ 
\text{with $\diam(U_j, \spa d_N^{T_1}) < \vep$ for all $ 1 \spa \leq j \spa \leq n$}
\end{array}
\right\}.
\end{align*}
(Letting $\Omega = X_1$, the above defines the standard topological pressure $P(f)$ on $(X_1, T_1)$. The topological entropy $h_{\mathrm{top}}(T_1)$ is the value of $P(f)$ when $f \equiv 0$.) Let $\Omega \subset X_{i+1}$. If $P^{\boldsymbol{a}}_i$ is already defined, let
\begin{flalign*}
& P^{\boldsymbol{a}}_{i+1}(\Omega, \spa f,\spa N, \spa \vep) &
\end{flalign*} \\[-35pt]
\begin{align*}
&= \inf \left\{ \spa \sum_{j=1}^n \Big( P^{\boldsymbol{a}}_i(\pi_i^{-1}(U_j), \spa f, \spa N, \spa \vep) \Big)^{a_i} \spa \middle|
\begin{array}{l}
\text{$n \in \mathbb{N}$, $\{U_j\}_{j=1}^n$ is an open cover of $\Omega$} \\ 
\text{with $\diam(U_j, \spa d_N^{T_{i+1}}) < \vep$ for all $ 1 \spa \leq j \spa \leq n$}
\end{array}
\right\}.
\end{align*} 
We define the \textbf{topological pressure of ${\boldsymbol{a}}$-exponent} $P^{\boldsymbol{a}}(f)$ by \\[1pt]
\[ P^{\boldsymbol{a}}(f) = \lim_{\vep \to 0} \left( \lim_{N \to \infty} \frac{\log{P^{\boldsymbol{a}}_r(X_r, \spa f, \spa N, \spa \vep)}}{N} \right). \] \\[1pt]
This limit exists since $\log{P^{\boldsymbol{a}}_r(X_r, \spa f, \spa N, \spa \vep)}$ is sub-additive in $N$ and non-decreasing as $\vep$ tends to $0$. When we want to clarify the maps $T_i$ and $\pi_i$ used in the definition of $P^{\boldsymbol{a}}(f)$, we will denote it by $P^{\boldsymbol{a}}(f, \spa \boldsymbol{T})$ or $P^{\boldsymbol{a}}(f, \spa \boldsymbol{T}, \spa \boldsymbol{\pi})$ with $\boldsymbol{T}=(T_i)_{i=1}^r$ and $\boldsymbol{\pi} = (\pi_i)_{i=1}^r$.

From ${\boldsymbol{a}}=(a_1, \spa a_2, \spa \cdots, a_{r-1})$, we define a probability vector (i.e., all entries are non-negative, and their sum is 1) $\boldsymbol{w_a} = (w_1, \spa\cdots, \spa w_r)$ by
\begin{eqnarray} \label{definition: w}
\left\{
\begin{array}{l}
w_1 = a_1 a_2 a_3 \cdots a_{r-1}\\
w_2 = (1-a_1) a_2 a_3 \cdots a_{r-1} \\
w_3 = (1-a_2) a_3 \cdots a_{r-1} \\
\hspace{50pt} \vdots \\
w_{r-1} = (1-a_{r-2}) a_{r-1} \\
w_r = 1- a_{r-1}
\end{array}
\right..
\end{eqnarray}
Let
\begin{gather*}
\pi^{(0)} = \mathrm{id}_{X_1}: X_1 \to X_1, \\
\pi^{(i)} = \pi_i \circ \pi_{i-1} \circ \cdots \circ \pi_1: X_1 \to X_{i+1}.
\end{gather*}
We can now state the main result of this paper.
\begin{theorem} \label{theorem: main theorem}
Let $(X_i, \spa T_i)$ ($i=1, \spa 2, \spa \ldots, \spa r$) be dynamical systems and $\pi_i: X_i \rightarrow X_{i+1} \hspace{5pt} (i=1, \spa 2, \spa ... \spa , \spa r-1)$ factor maps. For any continuous function $f: X_1 \to \mathbb{R}$,
\begin{equation} \label{equation: main result}
P^{\boldsymbol{a}}(f) = \sup_{\mu \in \mathscr{M}^{T_1}(X_1)} \left( \sum_{i=1}^r w_i h_{{\pi^{(i-1)}}_*\mu} (T_i) + w_1 \int_{X_1} f d\mu \right).
\end{equation}
\end{theorem}

We define $P^{\boldsymbol{a}}_\mathrm{var}(f)$ to be the right-hand side of this equation. Then we need to prove
\[ P^{\boldsymbol{a}}(f) = P^{\boldsymbol{a}}_\mathrm{var}(f). \]

\section{Preparation}
\subsection{Basic properties and tools} \label{subsection: basic properties and tools}
Let $(X_i, \spa T_i)$ ($i=1, \spa 2, \spa \ldots, \spa r$) be dynamical systems, $\pi_i: X_i \rightarrow X_{i+1} \hspace{5pt} (i=1, \spa 2, \spa \ldots \spa , \spa r-1)$ factor maps, $\boldsymbol{a} = (a_1, \cdots, a_{r-1}) \in [0, 1]^{r-1}$, and $f: X_1 \to \mathbb{R}$ a continuous function.

We will use the following notions in sections \ref{subsection: zero-dimensional principal extension} and \ref{section: proof of Pvar is smaller}.
\begin{definition} \label{definition: chain of covers}
Consider a cover $\flo{F}^{(i)}$ of $X_i$ for each $i$. For a natural number $N$ and a positive number $\vep$, the family $(\flo{F}^{(i)})_i$ is said to be \textbf{a chain of ($\boldsymbol{N}$, $\boldsymbol{\vep}$)-covers} of $(X_i)_i$ if the following conditions are true:
\begin{enumerate}
\item[(1)] For every $i$ and $V \in \flo{F}^{(i)}$, we have $\diam(V, d^{(i)}_N) < \vep$.
\item[(2)] For each $1 \leq i \leq r-1$ and $U \in \flo{F}^{(i+1)}$, there is $\flo{F}^{(i)}(U) \subset \flo{F}^{(i)}$ such that
\[ \pi_i^{-1}(U) \subset \bigcup \flo{F}^{(i)}(U) \]
and
\[ \flo{F}^{(i)} = \bigcup_{U \in \flo{F}^{(i+1)}} \flo{F}^{(i)}(U). \]
\end{enumerate}
Moreover, if all the elements of each $\flo{F}^{(i)}$ are open/closed/compact, we call $(\flo{F}^{(i)})_i$ \textbf{a chain of open/closed/compact ($\boldsymbol{N}$, $\boldsymbol{\vep}$)-covers} of $(X_i)_i$.
\end{definition}

\begin{remark} \label{remark: chains of covers}
Note that we can rewrite $P^{\boldsymbol{a}}_r(X_r, \spa f, \spa N, \spa \vep)$ using chains of open covers as follows. For a chain of ($N$, $\vep$)-covers $(\flo{F}^{(i)})_i$ of $(X_i)_i$, let \\[-10pt]
\begin{flalign*}
& \flo{P}^{\boldsymbol{a}}\left( f, \spa N, \spa \vep, \spa (\flo{F}^{(i)})_i \right) &
\end{flalign*} \\[-35pt]
\begin{align*}
&= \sum_{U^{(r)} \in \flo{F}^{(r)}} \left( \sum_{U^{(r-1)} \in \flo{F}^{(r-1)}(U^{(r)})} \left( \cdots \left( \sum_{U^{(1)} \in \flo{F}^{(1)}(U^{(2)})} e^{\sup_{U^{(1)}}S_Nf} \right)^{a_1} \cdots \right)^{a_{r-2}}\right)^{a_{r-1}}.
\end{align*}
Then
\begin{flalign*}
& P^{\boldsymbol{a}}_r(X_r, \spa f, \spa N, \spa \vep) &
\end{flalign*}
\begin{align*}
&= \inf{ \left\{ \flo{P}^{\boldsymbol{a}}\left( f, \spa N, \spa \vep, \spa (\flo{F}^{(i)})_i \right) \spa \middle| \spa (\flo{F}^{(i)})_i \text{ is a chain of open ($N$, $\vep$)-covers of $(X_i)_i$ } \right\} }.
\end{align*}
\end{remark}

Just like the classic notion of pressure, we have the following property.
\begin{lemma} \label{lemma: multiplication penetrates}
For any natural number m,
\[ P^{\boldsymbol{a}}(S_m^{T_1}f, \spa \boldsymbol{T}^m) = mP^{\boldsymbol{a}}(f, \spa \boldsymbol{T}), \]
{\it where} $\boldsymbol{T}^m = ({T_i}^m)_{i=1}^r$.
\end{lemma}
\begin{proof}
Fix $\vep > 0$.
It is obvious from the definition of $P^{\boldsymbol{a}}_1$ that for any 
$\Omega_1 \subset X_1$ and a natural number N,
\[ P^{\boldsymbol{a}}_1(\Omega_1, \spa S_m^{T_1}f, \spa \boldsymbol{T}^m, \spa N, \spa \vep)
\leq P^{\boldsymbol{a}}_1(\Omega_1, \spa f, \spa \boldsymbol{T}, \spa mN, \spa \vep). \]
Let $\Omega_{i+1} \subset X_{i+1}$. By induction on $i$, we have
\[ P^{\boldsymbol{a}}_i(\Omega_{i+1}, \spa S_m^{T_1}f, \spa \boldsymbol{T}^m\!, \spa N, \spa \vep)
\leq P^{\boldsymbol{a}}_i(\Omega_{i+1}, \spa f, \spa \boldsymbol{T}, \spa mN, \spa \vep). \]
Thus,
\begin{equation} \label{inequality: one side for multiplication penetrates}
P^{\boldsymbol{a}}_r(S_m^{T_1}f, \spa \boldsymbol{T}^m\!, \spa N, \spa \vep)
\leq P^{\boldsymbol{a}}_r(f, \spa \boldsymbol{T}, \spa mN, \spa \vep).
\end{equation}

There exists $0 < \delta < \vep$ such that for any $1 \leq i \leq r$,
\[ d^{(i)}(x, \spa y) < \delta \implies d_m^{T_i}(x, y) < \vep \qquad (\spa {\it for} \hspace{4pt}x, y \in X_i). \]
Then
\begin{equation} \label{equation: distance compared}
d_N^{T_i^m}(x, \spa y) < \delta \implies d_{mN}^{T_i}(x, y) < \vep \hspace{10pt}(\spa {\it for} \hspace{4pt} x, \spa y \in X_i \hspace{4pt} \it{and} \hspace{4pt} 1 \leq i \leq r).
\end{equation}
Let $i=1$ in (\refeq{equation: distance compared}), then we have for any $\Omega_1 \subset X_1$,
\[ P^{\boldsymbol{a}}_1(\Omega_1, \spa f, \spa \boldsymbol{T}, \spa mN, \spa \vep)
\leq P^{\boldsymbol{a}}_1(\Omega_1, \spa S_m^{T_1}f, \spa \boldsymbol{T}^m\!, \spa N, \spa \delta). \]
Take $\Omega_{i+1} \subset X_{i+1}$. Again by induction on $i$ and by (\refeq{equation: distance compared}), we have
\[ P^{\boldsymbol{a}}_i(\Omega_{i+1}, \spa f, \spa \boldsymbol{T}, \spa mN, \spa \vep)
\leq P^{\boldsymbol{a}}_i(\Omega_{i+1}, \spa S_m^{T_1}f, \spa \boldsymbol{T}^m\!, \spa N, \spa \delta). \]
Hence,
\[ P^{\boldsymbol{a}}_r(f, \spa \boldsymbol{T}, \spa mN, \spa \vep)
\leq P^{\boldsymbol{a}}_r(S_m^{T_1}f, \spa \boldsymbol{T}^m\!, \spa N, \spa \delta). \]

Combining with (\refeq{inequality: one side for multiplication penetrates}) we have
\[ P^{\boldsymbol{a}}_r(S_m^{T_1}f, \spa \boldsymbol{T}^m\!, \spa N, \spa \vep)
\leq P^{\boldsymbol{a}}_r(f, \spa \boldsymbol{T}, \spa mN, \spa \vep)
\leq P^{\boldsymbol{a}}_r(S_m^{T_1}f, \spa \boldsymbol{T}^m\!, \spa N, \spa \delta). \]
Therefore,
\[ P^{\boldsymbol{a}}(S^{T_1}_m f, \spa \boldsymbol{T}^m) = mP^{\boldsymbol{a}}(f, \spa \boldsymbol{T}). \]
\end{proof}

We will later use the following standard lemma of calculus.
\begin{lemma} \label{lemma: calculus}
\begin{minipage}[t]{\linewidth-\widthof{\the\csname thm@headfont\endcsname Theorem \thetheorem. }}
\begin{enumerate}[leftmargin=*]
	\item For $0 \leq a \leq 1$ and non-negative numbers $x, y$,
		\[ (x+y)^a \leq x^a+y^a. \]
	\item Suppose that non-negative real numbers $p_1, p_2, \ldots, p_n$ satisfy $\sum_{i=\mathrm{1}}^n p_i = \mathrm{1}$.
		Then for any real numbers $x_1, x_2, \ldots, x_n$ we have
		\[ \sum_{i=1}^n \left( -p_i \log{p_i} +x_i p_i \right) \leq \log{\sum_{i=1}^n e^{x_i}}. \]
		In particular, letting $x_1=x_2=\cdots=x_n=0$ gives
		\[ \sum_{i=1}^n(-p_i \log{p_i}) \leq \log{n}. \]
		Here, $0 \cdot \log{0}$ is defined as $0$.
\end{enumerate}%
\end{minipage}
\end{lemma}
The proof for (1) is elementary. See \cite[\S9.3, Lemma 9.9]{Walters} for (2).

\subsection{Measure theoretic entropy} \label{subsection: measure theoretic entropy}
In this subsection, we will introduce the classical measure-theoretic entropy (a.k.a. Kolmogorov-Sinai entropy) and state some of the basic lemmas we need to prove Theorem \ref{theorem: main theorem}. The main reference is the book of Walters \cite{Walters}.

Let $(X, T)$ be a dynamical system and $\mu \in \flo{M}^T(X)$. A set $\flo{A} = \{A_1, \ldots, A_n\}$ is called a finite partition of X with measurable elements if $X = A_1 \cup \dots \cup A_n$, each $A_i$ is a measurable set, and $A_i \cap A_j = \varnothing$ for $i \ne j$. In this paper, a partition is always finite and consists of measurable elements.

Let $\flo{A}$ and $\flo{A}'$ be partitions of $X$. We define a new partition $\flo{A} \vee \flo{A}'$ by
\[ \flo{A} \vee \flo{A}' = \left\{ A \cap A' \spa | \spa A \in \flo{A} \text{ and } A' \in \flo{A}' \right\}. \]
For a natural number $N$, we define a refined partition $\flo{A}_N$ of $\flo{A}$ by
\[ \flo{A}_N = \flo{A} \vee T^{-1}\flo{A} \vee T^{-2}\flo{A} \vee \cdots \vee T^{-(N-1)}\flo{A}, \] 
where $T^{-i}\flo{A} = \left\{ T^{-i}(A) \spa | \spa A \in \flo{A} \right\}$ is a partition for $i \in \mathbb{N}$.

For a partition $\flo{A}$ of $X$, let
\[ H_\mu(\flo{A}) = - \sum_{A \in \flo{A}} \mu(A) \log{(\mu(A))}. \]
We set
\[ h_\mu(T, \flo{A}) = \lim_{N \to \infty} \frac{H_\mu(\flo{A}_N)}{N}. \]
This limit exists since $H_\mu(\flo{A}_N)$ is sub-additive in $N$. The \textbf{measure theoretic entropy} $h_\mu(T)$ is defined by
\[ h_\mu(T) = \sup \left\{ h_\mu(T, \flo{A}) \spa | \spa \flo{A} \text{ is a partition of } X \right\}. \]

Let $\flo{A}$ and $\flo{A}'$ be partitions. Their \textbf{conditional entropy} is defined by
\[ H_\mu(\flo{A} | \flo{A}') = - \sum_{\substack{A' \in \flo{A}' \\ \mu(A') \ne 0}} \mu(A') \sum_{A \in \flo{A}} \frac{\mu(A \cap A')}{\mu(A')} \log{\left( \frac{\mu(A \cap A')}{\mu(A')} \right)}. \]

\begin{lemma} \label{lemma: properties of entropy}
\begin{minipage}[t]{\linewidth-\widthof{\the\csname thm@headfont\endcsname Theorem \thetheorem. }}
\begin{enumerate}[leftmargin=*]
\item $H_\mu(\flo{A})$ is sub-additive in $\flo{A}$: i.e., for partitions $\flo{A}$ and $\flo{A}'$,
\[ H_\mu(\flo{A} \vee \flo{A}') \leq H_\mu(\flo{A}) + H_\mu(\flo{A}'). \]
\item $H_\mu(\flo{A})$ is concave in $\mu$: i.e., for $\mu, \nu \in \flo{M}^T(X)$ and $0 \leq t \leq 1$,
\[ H_{(1-t)\mu+t\nu}(\flo{A}) \geq (1-t)H_\mu(\flo{A}) + tH_\nu(\flo{A}). \]
\item For partitions $\flo{A}$ and $\flo{A}'$,
		\[ h_\mu(T, \flo{A}) \leq h_\mu(T, \flo{A}') + H_\mu(\flo{A}' | \flo{A}). \]
\end{enumerate}%
\end{minipage}
\end{lemma}
For the proof confer \cite[Theorem 4.3 (viii), \S4.5]{Walters} for (1), \cite[Remark \S8.1]{Walters} for (2), and \cite[Theorem 4.12, \S4.5]{Walters} for (3).

\subsection{Zero-dimensional principal extension} \label{subsection: zero-dimensional principal extension}
Here we will see how we can reduce the proof of $P^{\boldsymbol{a}}(f) \leq P^{\boldsymbol{a}}_\mathrm{var}(f)$ to the case where all dynamical systems are zero-dimensional.

First, we review the definitions and properties of (zero-dimensional) principal extension. The introduction here closely follows Tsukamoto's paper \cite{Tsukamoto} and references the book of Downarowicz \cite{Downarowicz}. Suppose $\pi: (Y, \spa S) \rightarrow (X, \spa T)$ is a factor map between dynamical systems. Let $d$ be a metric on $Y$. We define the \textbf{conditonal topological entropy} of $\pi$ by
\[ h_\mathrm{top}(Y, S \spa | \spa X, T) = \lim_{\vep \to 0} \left( \lim_{N \to \infty} \frac{\sup_{x \in X} \log{\#(\pi^{-1}(x), N, \vep)}}{N} \right). \]
Here,
\begin{align*}
\#(\pi^{-1}(x), \spa N, \spa \vep)
&= \min \left\{ \spa n \in \mathbb{N} \spa \middle|
\begin{array}{l}
\text{There exists an open cover $\{U_j\}_{j=1}^n$ of $\pi^{-1}(x)$} \\ 
\text{with $\diam(U_j, \spa d_N) < \vep$ for all $ 1 \spa \leq j \spa \leq n$}
\end{array}
\right\}.
\end{align*}

A factor map $\pi: (Y, \spa S) \rightarrow (X, \spa T)$ between dynamical systems is said to be a \textbf{principal factor map} if
\[ h_\mathrm{top}(Y, S \spa | \spa X, T) = 0. \]
Also, $(Y, \spa S)$ is called a \textbf{principal extension} of $(X, \spa T)$.

The following theorem is from \cite[Corollary 6.8.9]{Downarowicz}.
\begin{theorem}
Suppose $\pi: (Y, \spa S) \rightarrow (X, \spa T)$ is a principal factor map. Then $\pi$ preserves measure-theoretic entropy, namely,
\[ h_\mu(S) = h_{\pi_*\mu}(T) \]
for any $S$-invariant probability measure $\mu$ on Y.
\end{theorem}
More precisely, it is proved in \cite[Corollary 6.8.9]{Downarowicz} that $\pi$ is a principal factor map if and only if it preserves measure-theoretic entropy.

Suppose $\pi: (X_1, T_1) \rightarrow (X_2, T_2)$ and $\phi: (Y, S) \rightarrow (X_2, T_2)$ are factor maps between dynamical systems. We define a fiber product $(X_1 \times_{X_2}^{} Y, \spa T_1 \times S)$ of $(X_1, T_1)$ and $(Y, S)$ over $(X_2, T_2)$ by
\[ X_1 \times_{X_2}^{} Y = \left\{ (x, y) \in X_1 \times Y \spa \middle| \spa \pi(x) = \phi(y) \right\}, \]
\[ T_1 \times S: X_1 \times_{X_2}^{} Y \ni (x, y) \longmapsto \left(T_1(x), S(y)\right) \in X_1 \times_{X_2}^{} Y. \]
We have the following commutative diagram:
\begin{equation} \label{diagram: row extension one}
\begin{gathered}
 \xymatrix@C=40pt@R=36pt{
	X_1 \times_{X_2}^{} Y \ar[d]_-{\pi'} \ar[r]^-\psi & X_1 \ar[d]^-\pi \\
	Y \ar[r]_-\phi & X_2 }
\end{gathered}
\end{equation}
Here $\pi'$ and $\psi$ are restrictions of the projections onto $Y$ and $X_1$, respectively:
\[ \pi':  X_1 \times_{X_2}^{} Y \ni (x, y) \longmapsto y \in Y, \]
\[ \psi: X_1 \times_{X_2}^{} Y \ni (x, y) \longmapsto x \in X_1. \]
Since $\pi$ and $\phi$ are surjective, both $\pi'$ and $\psi$ are factor maps.

\begin{lemma} \label{lemma: square extension}
If $\phi$ is a principal extension in the diagram (\refeq{diagram: row extension one}), then $\psi$ is also a principal extension.
\end{lemma}
\begin{proof}
Let $d^1$ and $d^Y$ be metrics on $X_1$ and $Y$, respectively. Define a metric $\widetilde{d}$ on $X_1 \times_{X_2}^{} Y$ by
\[ \widetilde{d}\big( (x, y), (x', y') \big) = \max{\{ d^1(x, x'), d^Y(y, y') \}}. \]
Let $x \in X_1$. We have
\begin{equation*}
\psi^{-1}(x) = \{x\} \times \left\{ y \in Y \spa \middle| \spa \pi(x) = \phi(y) \right\} = \{x\} \times \phi^{-1}(\pi(x)),
\end{equation*}
which in turn implies $\widetilde{d}|_{\psi^{-1}(x)} = d^Y|_{\phi^{-1}(\pi(x))}$.
Then the metric space $(\psi^{-1}(x), \widetilde{d}_N)$ is isometric to $(\phi^{-1}(\pi(x)), d^Y_N)$ for any natural number $N$. Therefore for any $\vep > 0$,
\[ \#(\psi^{-1}(x), N, \vep) = \#(\phi^{-1}(\pi(x)), N, \vep). \]
Since $\pi$ is surjective,
\[ \sup_{x \in X_1} \#(\psi^{-1}(x), N, \vep) = \sup_{x \in X_1} \#(\phi^{-1}(\pi(x)), N, \vep) = \sup_{y \in Y} \#(\phi^{-1}(y), N, \vep). \]
Hence,
\[ h_\mathrm{top}(X_1 \times_{X_2}^{} Y, T_1 \times S \spa | \spa X_1, T_1) = h_\mathrm{top}(Y, S \spa | \spa X_2, T_2) = 0. \]
\end{proof}

A dynamical system $(Y, S)$ is said to be \textbf{zero-dimensional} if there is a clopen basis of the topology of $Y$, where clopen means any element in the basis is both closed and open. A basic example of a zero-dimensional dynamical system is the Cantor set $\{ 0, 1 \}^\mathbb{N}$ with the shift map.

A principal extension $(Y, \spa S)$ of $(X, \spa T)$ is called a \textbf{zero-dimensional principal extension} if $(Y, S)$ is zero-dimensional. The following important theorem can be found in \cite[Theorem 7.6.1]{Downarowicz}.
\begin{theorem} \label{theorem: zero dimensional principal extension}
For any dynamical system, there is a zero-dimensional principal extension.
\end{theorem}

Let $(Y_i, \spa R_i)$ ($i=1, \spa 2, \spa \ldots, \spa m$) be dynamical systems, $\pi_i: Y_i \rightarrow Y_{i+1} \hspace{5pt} (i=1, \spa 2, \spa ... \spa , \spa m-1)$ factor maps, and $\boldsymbol{a} = (a_1, \cdots, a_{m-1}) \in [0, 1]^{m-1}$. Fix $2 \leq k \leq m-1$ and take a zero-dimensional principal extension $\phi_k: (Z_k, S_k) \rightarrow (Y_k, R_k)$. 
For each $1 \leq i \leq k-1$, let $(Y_i \times_{Y_k} Z_k, R_i \times S_k)$ be the fiber product and $\phi_i: Y_i \times_{Y_k} Z_k \rightarrow Y_i$ be the restriction of the projection as in the earlier definition. We have
\[ \xymatrix@C=60pt{
Y_i \times_{Y_k} Z_k \ar[r]^-{\phi_i} \ar[d] & Y_i \ar[d]^-{\pi_{k-1} \circ \pi_{k-2} \circ \cdots \circ \pi_i} \\
Z_k \ar[r]_-{\phi_k} & Y_k } \]
By Lemma \ref{lemma: square extension}, $\phi_i$ is a principal factor map.
We define $\Pi_i: Y_i \times_{Y_k} Z_k \rightarrow Y_{i+1} \times_{Y_k} Z_k$ by $\Pi_i(x, y) = \left( \pi_i(x), y \right)$ for each $i$. Then we have the following commutative diagram:
\begin{equation} \label{diagram: row extension}
\begin{gathered}
\xymatrix@C=60pt@R=18pt{
Y_1 \times_{Y_k}^{} Z_k \ar[r]^-{\phi_1} \ar[d]_-{\Pi_1} & Y_1 \ar[d]^-{\pi_1} \\
Y_2 \times_{Y_k}^{} Z_k \ar[r]^-{\phi_2} \ar[d]_-{\Pi_2} & Y_2 \ar[d]^-{\pi_2} \\
\vdots \ar[d]_-{\Pi_{k-2}} &  \vdots \ar[d]^-{\pi_{k-2}} \\
Y_{k-1} \times_{Y_k}^{} Z_k \ar[r]^-{\phi_{k-1}} \ar[d]_-{\Pi_{k-1}} & Y_{k-1} \ar[d]^-{\pi_{k-1}} \\
Z_k \ar[rd]_-{\pi_k \circ \phi_k} \ar[r]^-{\phi_k} & Y_k \ar[d]^-{\pi_k} \\
 	& Y_{k+1} \ar[d]^-{\pi_{k+1}} \\
	& \vdots \ar[d]^-{\pi_{m-1}} \\
	& Y_m }
\end{gathered}
\end{equation}
Let
\begin{equation*}
\begin{gathered}
(Z_i, S_i) = (Y_i \times^{}_{Y_k} Z_k, R_i \times S_k) \text{ for $1 \leq i \leq k-1$},
\hspace{6pt} (Z_i, S_i) = (Y_i, R_i) \text{   for $k+1 \leq i \leq m$}, \\
\Pi_k = \pi_k \circ \phi_k: Z_k \rightarrow Y_{k+1}, \hspace{6pt}
\Pi_i = \pi_i: Z_i \rightarrow Z_{i+1} \text{  for $k+1 \leq i \leq m-1$}, \\
\phi_i = \mathrm{id}_{Z_i}: Z_i \rightarrow Z_i  \text{   for $k+1 \leq i \leq m$}.
\end{gathered}
\end{equation*}

\begin{lemma} \label{lemma: pressure inequality}
In the settings above,
\[ P^{\boldsymbol{a}}_{\mathrm{var}}(f, \boldsymbol{R}, \boldsymbol{\pi}) \geq P^{\boldsymbol{a}}_{\mathrm{var}}(f \circ \phi_1, \boldsymbol{S}, \boldsymbol{\Pi}) \]
and
\[ P^{\boldsymbol{a}}(f, \boldsymbol{R}, \boldsymbol{\pi}) \leq P^{\boldsymbol{a}}(f \circ \phi_1, \boldsymbol{S}, \boldsymbol{\Pi}). \]
Here, $\boldsymbol{R} = (R_i)_i$, $\boldsymbol{\pi} = (\pi_i)_i$, $\boldsymbol{S} = (S_i)_i$ and $\boldsymbol{\Pi} = (\Pi_i)_i$.
\end{lemma}
\begin{proof}
We remark that the following proof does not require $Z_k$ to be zero-dimensional. Let
\begin{gather*}
\pi^{(0)} = \mathrm{id}_{Y_1}: Y_1 \to Y_1, \\
\pi^{(i)} = \pi_i \circ \pi_{i-1} \circ \cdots \circ \pi_1: Y_1 \to Y_{i+1}
\end{gather*}
and
\begin{gather*}
\Pi^{(0)} = \mathrm{id}_{Z_1}: Z_1 \to Z_1, \\
\Pi^{(i)} = \Pi_i \circ \Pi_{i-1} \circ \cdots \circ \Pi_1: Z_1 \to Z_{i+1}.
\end{gather*}
Let $\nu \in \flo{M}^{S_1}(Y_1)$ and $1 \leq i \leq m$. Since all the horizontal maps in \eqref{diagram: row extension} are principal factor maps, we have
\[ h_{{\Pi^{(i-1)}}_*\nu}(S_i) = h_{(\phi_i)_*{\Pi^{(i-1)}}_*\nu}(R_i) = h_{{\pi^{(i-1)}}_*(\phi_1)_*\nu}(R_i). \]
It follows that
\begin{align*}
P^{\boldsymbol{a}}_{\mathrm{var}}(f \circ \phi_1, \boldsymbol{S}, \boldsymbol{\Pi})
&= \sup_{\nu \in \mathscr{M}^{S_1}(Z_1)} \left( \sum_{i=1}^{m} w_i h_{{\Pi^{(i-1)}}_*\nu} (S_i) + w_1 \int_{Z_1} f \circ \phi_1 d\nu \right) \\
&= \sup_{\nu \in \mathscr{M}^{S_1}(Z_1)} \left( \sum_{i=1}^{m} w_i h_{{\pi^{(i-1)}}_*(\phi_1)_*\nu}(R_i) + w_1 \int_{Y_1} f d\big( (\phi_1)_*\nu \big) \right) \\
&\leq \sup_{\mu \in \mathscr{M}^{T_1}(Y_1)} \left( \sum_{i=1}^{m} w_i h_{{\pi^{(i-1)}}_*\mu} (R_i) + w_1 \int_{Y_1} f d\mu \right) \\[4pt]
&= P^{\boldsymbol{a}}_{\mathrm{var}}(f, \boldsymbol{R}, \boldsymbol{\pi}).
\end{align*}
(The reversed inequality is generally true by the surjectivity of factor maps, yielding equality. However, we do not use this fact.)

Let $d^i$ be a metric on $Y_i$ for each $i$ and $\widetilde{d^k}$ a metric on $Z_k$. We define a metric $\widetilde{d^i}$ on $(Z_i, S_i)$ for $1 \leq i \leq k-1$ by
\[ \widetilde{d^i}\big( (x_1, y_1), (x_2, y_2) \big) = \max{\{ d^i(x_1, x_2), \widetilde{d^k}(y_1, y_2) \}} \text{ \hspace{4pt} $\big($$(x_1, y_1), (x_2, y_2) \in Z_i = Y_i \times_{Y_k} Z_k$$\big)$ }. \]
Set $\widetilde{d^i} = d^i$ for $k+1 \leq i \leq m$. Take an arbitrary positive number $\vep$. There exists $0 < \delta < \vep$ such that for every $1 \leq i \leq m$,
\begin{equation} \label{inequality: distance of factors compared}
\widetilde{d^i}(x, y) < \delta \implies d^i( \phi_i(x), \phi_i(y) ) < \vep \quad (x, y \in Z_i).
\end{equation}

Let $N$ be a natural number. We claim that
\[ P^{\boldsymbol{a}}_r(f, \boldsymbol{R}, \boldsymbol{\pi}, N, \vep)
\leq P^{\boldsymbol{a}}_r(f \circ \phi_1, \boldsymbol{S}, \boldsymbol{\Pi}, N, \delta). \]
Take $M > 0$ with
\[ P^{\boldsymbol{a}}_r(f \circ \phi_1, \boldsymbol{S}, \boldsymbol{\Pi}, N, \delta) < M. \]
Then there exists a chain of open ($N$, $\delta$)-covers $(\flo{F}^{(i)})_i$ of $(Z_i)_i$ (see Definition \ref{definition: chain of covers} and Remark \ref{remark: chains of covers}) with
\[ \flo{P}^{\boldsymbol{a}}\left(f \circ \phi_1, \spa \boldsymbol{S}, \spa \boldsymbol{\Pi}, \spa N, \spa \delta, \spa (\flo{F}^{(i)})_i \right) < M. \]
We can find a compact set $C_U \subset U$ for each $U \in \flo{F}^{(m)}$ such that $\bigcup_{U \in \flo{F}^{(m)}} C_U = Z_m$. Let $\flo{K}^{(m)} := \{ C_U \spa | \spa U \in \flo{F}^{(m)} \}$. Since $\Pi_{m-1}^{-1}(C_U) \subset \Pi_{m-1}^{-1}(U)$ is compact for each $U \in \flo{F}^{(m)}$, we can find a compact set $E_V \subset V$ for each $V \in \flo{F}^{(m-1)}(U)$ such that $\Pi_{m-1}^{-1}(C_U) \subset \bigcup_{V \in \flo{F}^{(k)}(U)} E_V$. Let $\flo{K}^{(m-1)}(C_U) := \{ E_V \spa | \spa V \in \flo{F}^{(m-1)}(U) \}$ and $\flo{K}^{(m-1)} := \bigcup_{C \in \flo{K}^{(m)}} \flo{K}^{(m-1)}(C)$. We continue likewise and obtain a chain of compact ($N$, $\delta$)-covers $(\flo{K}^{(i)})_i$ of $(Z_i)_i$ with
\[ \flo{P}^{\boldsymbol{a}}\left(f \circ \phi_1, \spa \boldsymbol{S}, \spa \boldsymbol{\Pi}, \spa N, \spa \delta, \spa (\flo{K}^{(i)})_i \right)
\leq \flo{P}^{\boldsymbol{a}}\left(f \circ \phi_1, \spa \boldsymbol{S}, \spa \boldsymbol{\Pi}, \spa N, \spa \delta, \spa (\flo{F}^{(i)})_i \right) < M. \]
Let $\phi_i(\flo{K}^{(i)}) = \left\{ \phi_i(C) \spa \middle| \spa C \in \flo{K}^{(i)} \right\}$ for each $i$. Note that for any $\Omega \subset Z_i$,
\[ \pi_{i-1}^{-1} ( \phi_{i} ( \Omega ) ) = \phi_{i-1} ( \Pi_{i-1}^{-1} ( \Omega ) ). \]
This and (\refeq{inequality: distance of factors compared}) assure that $(\phi_i(\flo{K}^{(i)}))_i$ is a chain of compact ($N$, $\vep$)-covers of $(Y_i)_i$. We have
\begin{align*}
\flo{P}^{\boldsymbol{a}}\left( f, \spa \boldsymbol{R}, \spa \boldsymbol{\pi}, \spa N, \spa \vep, \spa (\phi_i(\flo{K}^{(i)}))_i \right)
&=
\flo{P}^{\boldsymbol{a}}\left( f \circ \phi_1, \spa \boldsymbol{S}, \spa \boldsymbol{\Pi}, \spa N, \spa \delta, \spa (\flo{K}^{(i)})_i \right) < M.
\end{align*}

Since $f$ is continuous and each $\phi_i(\flo{K}^{(i)})$ is a closed cover, we can slightly enlarge each set in $\phi_i(\flo{K}^{(i)})$ and create a chain of open ($N$, $\vep$)-covers $(\flo{O}^{(i)})_i$ of $(Y_i)_i$ satisfying
\begin{align*}
\flo{P}^{\boldsymbol{a}}\left( f, \spa \boldsymbol{R}, \spa \boldsymbol{\pi}, \spa N, \spa \vep, \spa (\flo{O}^{(i)})_i \right) < M.
\end{align*}
Therefore,
\[ P^{\boldsymbol{a}}_r(f, \boldsymbol{R}, \boldsymbol{\pi}, N, \vep) \leq \flo{P}^{\boldsymbol{a}}\left( f, \spa \boldsymbol{R}, \spa \boldsymbol{\pi}, \spa N, \spa \vep, \spa (\flo{O}^{(i)})_i \right) < M. \]
Since $M > P^{\boldsymbol{a}}_r(f \circ \phi_1, \boldsymbol{S}, \boldsymbol{\Pi}, N, \delta)$ was chosen arbitrarily, we have 
\[ P^{\boldsymbol{a}}_r(f, \boldsymbol{R}, \boldsymbol{\pi}, N, \vep) \leq P^{\boldsymbol{a}}_r(f \circ \phi_1, \boldsymbol{S}, \boldsymbol{\Pi}, N, \delta). \]
This implies
\[ P^{\boldsymbol{a}}(f, \boldsymbol{R}, \boldsymbol{\pi}) \leq P^{\boldsymbol{a}}(f \circ \phi_1, \boldsymbol{S}, \boldsymbol{\Pi}). \]
\end{proof}

The following proposition reduces the proof of $P^{\boldsymbol{a}}(f) \leq P^{\boldsymbol{a}}_\mathrm{var}(f)$ in the next section to the case where all dynamical systems are zero-dimensional.

\begin{proposition} \label{proposition: zero-dimensional trick}
For all dynamical systems $(X_i, \spa T_i)$ ($i=1, \spa 2, \spa \ldots, \spa r$) and factor maps $\pi_i: X_i \rightarrow X_{i+1} \hspace{5pt} (i=1, \spa 2, \spa ... \spa , \spa r-1)$, there are zero-dimensional dynamical systems $(Z_i, \spa S_i)$ ($i=1, \spa 2, \spa \ldots, \spa r$) and factor maps $\Pi_i: Z_i \rightarrow Z_{i+1} \hspace{5pt} (i=1, \spa 2, \spa ... \spa , \spa r-1)$ with the following property; for every continuous function $f: X_1 \rightarrow \mathbb{R}$ there exists a continuous function $g: Z_1 \rightarrow \mathbb{R}$ with
\[ P^{\boldsymbol{a}}_{\mathrm{var}}(f, \boldsymbol{T}, \boldsymbol{\pi}) \geq P^{\boldsymbol{a}}_{\mathrm{var}}(g, \boldsymbol{S}, \boldsymbol{\Pi}) \]
and
\[ P^{\boldsymbol{a}}(f, \boldsymbol{T}, \boldsymbol{\pi}) \leq P^{\boldsymbol{a}}(g, \boldsymbol{S}, \boldsymbol{\Pi}). \]
\end{proposition}

\begin{proof}
We will first construct zero-dimensional dynamical systems $(Z_i, \spa S_i)$ ($i=1, \spa 2, \spa \ldots, \spa r$) and factor maps $\Pi_i: Z_i \rightarrow Z_{i+1} \hspace{5pt} (i=1, \spa 2, \spa ... \spa , \spa r-1)$ alongside the following commutative diagram of dynamical systems and factor maps: 

\begin{equation} \label{diagram: Xi to Zi}
\xymatrix@C=9pt@R=6pt{
	Z_1 \ar[r]^-{\phi_r} \ar[rdd]_{\Pi_1} & \ar[dd] & & & \cdots \hspace{90pt} \ar[r]^-{\phi_2} & X_1 \times_{X_r}^{} Z_r \ar[r]^-{\phi_1} \ar[dd]^{\pi_1^{(2)}} & X_1 \ar[dd]^{\pi_1} \\
	&&&&&& \\
	      & Z_2 \ar[r] \ar[rdd]_{\Pi_2} & & & \cdots \hspace{90pt} \ar[r]  & X_2 \times_{X_r}^{} Z_r \ar[r] \ar[dd]^{\pi_2^{(2)}} & X_2 \ar[dd]^{\pi_2} \\
	&&&&&& \\ 
	& & \ddots \hspace{4pt} \ar[rd]_{\Pi_{r-3}} & \ar[d]^{\pi^{(4)}_{r-3}} & \vdots \ar[d]^{\pi_{r-3}^{(3)}} & \vdots \ar[d]^{\pi_{r-3}^{(2)}} & \vdots \ar[d]^{\pi_{r-3}} \\ 
	      & & &	Z_{r-2} \ar[r] \ar[rddd]_{\Pi_{r-2}} & \left(X_{r-2} \times_{X_r}^{} Z_r\right) \times_{(X_{r-1} \times_{X_r} Z_r)}^{} Z_{r-1} \ar[r] \ar[ddd]^{\pi_{r-2}^{(3)}} & X_{r-2} \times_{X_r}^{} Z_r \ar[r] \ar[ddd]^{\pi_{r-2}^{(2)}} & X_{r-2} \ar[ddd]^{\pi_{r-2}} \\
	&&&&&& \\ 	&&&&&& \\
	      & & & & 	    Z_{r-1} \ar[r]^-{\psi_{r-1}} \ar[rdddd]_{\Pi_{r-1}} & X_{r-1} \times_{X_r}^{} Z_r \ar[r] \ar[dddd]^{\pi_{r-1}^{(2)}} & X_{r-1} \ar[dddd]^{\pi_{r-1}} \\
	&&&&&& \\	&&&&&& \\ 	&&&&&& \\
	& & & & & Z_r \ar[r]^-{\psi_r} \ar[rdd] & \hspace{2pt} X_r \ar[dd] \\
	&&&&&& \\
	& & & & & & \{*\} }
\end{equation}
where all the horizontal maps are principal factor maps.

By Theorem \ref{theorem: zero dimensional principal extension}, there is a zero-dimensional principal extension $\psi_r: (Z_r, S_r) \rightarrow (X_r, T_r)$. The set $\{*\}$ is the trivial dynamical system, and the maps $X_r \rightarrow \{*\}$ and $Z_r \rightarrow \{*\}$ send every element to $*$. For each $1 \leq i \leq r-1$, the map $X_i \times_{X_r}^{} Z_r \rightarrow X_i$ in the following diagram is a principal factor map by Lemma \ref{lemma: square extension}.
\[ \xymatrix@C=40pt@R=36pt{
	X_i \times_{X_r}^{} Z_r \ar[d] \ar[r] & X_i \ar[d]^-{\pi_{r-1} \circ \pi_{r-2} \circ \cdots \circ \pi_i} \\
	Z_r \ar[r]_-{\psi_r} & X_r } \]
For $1 \leq i \leq r-2$, define $\pi_i^{(2)}: X_i \times_{X_r}^{} Z_r \rightarrow X_{i+1} \times_{X_r}^{} Z_r$ by
\[ \pi_i^{(2)}(x, z) = (\pi_i(x), y). \]
Then every horizontal map in the right two rows of (\refeq{diagram: Xi to Zi}) is a principal factor map. Next, take a zero-dimensional principal extension $\psi_{r-1}: (Z_{r-1}, S_{r-1}) \rightarrow (X_{r-1} \times_{X_r}^{} Z_r, T_{r-1} \times S_r)$ and let $\Pi_{r-1} = \pi_{r-1}^{(2)} \circ \psi_{r-1}$. The rest of (\refeq{diagram: Xi to Zi}) is constructed similarly, and by Lemma \ref{lemma: square extension}, each horizontal map is a principal factor map.

Let $f: X_1 \rightarrow \mathbb{R}$ be a continuous map. Applying Lemma \ref{lemma: pressure inequality} to the right two rows of (\refeq{diagram: Xi to Zi}), we get
\[ P^{\boldsymbol{a}}_{\mathrm{var}}(f, \boldsymbol{T}, \boldsymbol{\pi}) \geq P^{\boldsymbol{a}}_{\mathrm{var}}(f \circ \phi_1, \boldsymbol{S^{(2)}}, \boldsymbol{\Pi^{(2)}}) \]
and
\[ P^{\boldsymbol{a}}(f, \boldsymbol{T}, \boldsymbol{\pi}) \leq P^{\boldsymbol{a}}(f \circ \phi_1, \boldsymbol{S^{(2)}}, \boldsymbol{\Pi^{(2)}}) \]
for $\boldsymbol{\Pi^{(2)}} = (\pi^{(2)}_i)_i$ and $\boldsymbol{S^{(2)}} = (T_i \times S_r)_i$. Again by Lemma \ref{lemma: pressure inequality},
\[ P^{\boldsymbol{a}}_{\mathrm{var}}(f \circ \phi_1, \boldsymbol{S^{(2)}}, \boldsymbol{\Pi^{(2)}}) \geq P^{\boldsymbol{a}}_{\mathrm{var}}(f \circ \phi_1 \circ \phi_2, \boldsymbol{S^{(3)}}, \boldsymbol{\Pi^{(3)}}) \]
and
\[ P^{\boldsymbol{a}}(f \circ \phi_1, \boldsymbol{S^{(2)}}, \boldsymbol{\Pi^{(2)}}) \leq P^{\boldsymbol{a}}(f \circ \phi_1 \circ \phi_2, \boldsymbol{S^{(3)}}, \boldsymbol{\Pi^{(3)}}) \]
where $\boldsymbol{\Pi^{(3)}} = \big( (\pi^{(3)}_i)_{i=1}^{r-2}, \Pi_{r-1} \big)$, and $\boldsymbol{S^{(3)}}$ is the collection of maps associated with $Z_r$ and the third row from the right of (\refeq{diagram: Xi to Zi}).
We continue inductively and obtain the desired inequalities, where $g$ is taken as $f \circ \phi_1 \circ \phi_2 \circ \cdots \circ \phi_r$.
\end{proof}

\section{\texorpdfstring{Proof of $P^{\boldsymbol{a}}(f) \leq P^{\boldsymbol{a}}_\mathrm{var}(f)$.}%
                               {subtitle}}
Let $\boldsymbol{a} = (a_1, \cdots, a_{r-1}) \in [0, 1]^{r-1}$. Recall that we defined $(w_1, \ldots, w_r)$ by
\begin{eqnarray*}
\left\{
\begin{array}{l}
w_1 = a_1 a_2 a_3 \cdots a_{r-1}\\
w_2 = (1-a_1) a_2 a_3 \cdots a_{r-1} \\
w_3 = (1-a_2) a_3 \cdots a_{r-1} \\
\hspace{50pt} \vdots \\
w_{r-1} = (1-a_{r-2}) a_{r-1} \\
w_r = 1- a_{r-1}
\end{array}
\right.
\end{eqnarray*}
and
$P^{\boldsymbol{a}}_\mathrm{var}(f)$ by
\[ P^{\boldsymbol{a}}_\mathrm{var}(f) = \sup_{\mu \in \mathscr{M}^{T_1}(X_1)} \left( \sum_{i=1}^r w_i h_{{\pi^{(i-1)}}_*\mu} (T_i) + w_1 \int_{X_1} f d\mu \right) \]
where
\begin{gather*}
\pi^{(0)} = \mathrm{id}_{X_1}: X_1 \to X_1, \\
\pi^{(i)} = \pi_i \circ \pi_{i-1} \circ \cdots \circ \pi_1: X_1 \to X_{i+1}.
\end{gather*}

The following theorem suffices by Theorem \ref{proposition: zero-dimensional trick} in proving $P^{\boldsymbol{a}}(f) \leq P^{\boldsymbol{a}}_\mathrm{var}(f)$ for arbitrary dynamical systems.
\begin{theorem} \label{theorem: first half of the main theorem}
Suppose $(X_i, \spa T_i)$ ($i=1, \spa 2, \spa \ldots, \spa r$) are zero-dimensional dynamical systems and $\pi_i: X_i \rightarrow X_{i+1} \hspace{5pt} (i=1, \spa 2, \spa ... \spa , \spa r-1)$ are factor maps. Then we have
\[ P^{\boldsymbol{a}}(f) \leq P^{\boldsymbol{a}}_\mathrm{var}(f) \]
for any continuous function $f: X_1 \rightarrow \mathbb{R}$.
\end{theorem}
\begin{proof}
Let $d^{(i)}$ be a metric on $X_i$ for each $i=1, 2, \ldots, r$. Take a positive number $\vep$ and a natural number $N$. First, we will backward inductively define a finite clopen partition $\flo{A}^{(i)}$ of $X_i$ for each $i$. Since $X_r$ is zero-dimensional, we can take a sufficiently fine finite clopen partition $\flo{A}^{(r)}$ of $X_r \spa$. That is, each $A \in \flo{A}^{(r)}$ is both open and closed, and $\diam(A, d^{(r)}_N) < \vep$. Suppose $\flo{A}^{(i+1)}$ is defined. For each $A \in \flo{A}^{(i+1)}$, take a clopen partition $\flo{B}(A)$ of $\pi_{i}^{-1} (A) \subset X_{i}$ such that any $B \in \flo{B}(A)$ satisfies $\diam(B, d^{(i)}_N) < \vep$. We let $\flo{A}^{(i)} = \bigcup_{A \in \flo{A}^{(i+1)}} \flo{B}(A)$. Then $\flo{A}^{(i)}$ is a finite clopen partition of $X_i$.
We define 
\[ \flo{A}^{(i)}_N = \flo{A}^{(i)} \vee T_i^{-1}\flo{A}^{(i)} \vee T_i^{-2}\flo{A}^{(i)} \vee \cdots
\vee T_i^{-(N-1)}\flo{A}^{(i)}. \]

We employ the following notations. For $i<j$ and $A \in \flo{A}^{(j)}_N$, let $\flo{A}^{(i)}_N(A)$ be the set of ``children'' of A;
\[ \flo{A}^{(i)}_N(A) = \left\{B \in \flo{A}^{(i)}_N \spa \middle| \spa \pi_{j-1} \circ \pi_{j-2} \circ \cdots \circ \pi_i(B) \subset A \right\}. \]
Also, for $B \in \flo{A}^{(i)}_N$ and $i<j$,  we denote by $\tpi_j B$ the unique ``parent'' of $B$ in $\flo{A}^{(j)}_N$;
\[ \tpi_j B = A \in \flo{A}^{(j)}_N \text{ such that } \pi_{j-1} \circ \pi_{j-2} \circ \cdots \circ \pi_{i}(B) \subset A. \]

We will evaluate $P^{\boldsymbol{a}}(f, \spa N, \spa \vep)$ from above using $\{ \flo{A}^{(i)} \}$. Let $A \in \flo{A}^{(2)}_N$, and start by setting
\[ Z^{(1)}_N (A) = \sum_{B \in \flo{A}^{(1)}_N(A)} e^{\sup_B S_N f}. \]
Let $A \in \flo{A}^{(i+1)}_N$. If $Z^{(i-1)}_N$ is already defined, set
\[ Z^{(i)}_N (A) = \sum_{B \in \flo{A}^{(i)}_N(A)} \left( Z^{(i-1)}_N (B) \right)^{a_{i-1}}. \]
We then define $Z_N$ by
\[ Z_N = \sum_{A \in \flo{A}^{(r)}_N} \left( Z^{(r-1)}_N (A) \right)^{a_{r-1}}. \]
It is straightforward from the construction that
\[ P^{\boldsymbol{a}}_r(X_r, \spa f, \spa N, \spa \vep) \leq Z_N. \]
Therefore, we only need to prove that there is a $T_1$-invariant probability measure $\mu$ on $X_1$ such that
\[ \sum_{i=1}^r w_i h_{{\pi^{(i-1)}}_* \mu} (T_i, \spa \flo{A}^{(i)}) + w_1 \int_{X_1} f d\mu \geq \lim_{N \to \infty} \frac{\log Z_N}{N}. \]

Since each $A \in \flo{A}^{(1)}_N$ is closed, we can choose a point $x_A \in A$ so that
\[ S_N f(x_A) = \sup_A S_N f. \]
We define a probability measure $\sigma^{}_N$ on $X_1$ by
\begin{equation*}
\sigma_N^{} =
\begin{multlined}[t][11.5cm]
\frac{1}{Z_N} \sum_{A \in \flo{A}^{(1)}_N} {Z_N^{(r-1)}(\tpi_{r}A)}^{a_{r-1}-1}{Z_N^{(r-2)}(\tpi_{r-1}A)}^{a_{r-2}-1} \\
\times \cdots \times {Z_N^{(2)}(\tpi_3 A)}^{a_{2}-1}{Z_N^{(1)}(\tpi_2 A)}^{a_{1}-1}e^{S_N f(x_A)} \delta_{x_A}
\end{multlined}
\end{equation*}
where $\delta_{x^{}_A}$ is the Dirac measure at $x_A$. This is indeed a probability measure on $X_1$ since
\begin{align*}
\sigma_N^{}(X_1) &=
\begin{multlined}[t][11.5cm]
\frac{1}{Z_N} \sum_{A \in \flo{A}^{(1)}_N} {Z_N^{(r-1)}(\tpi_{r}A)}^{a_{r-1}-1}{Z_N^{(r-2)}(\tpi_{r-1}A)}^{a_{r-2}-1} \\
\times \cdots \times {Z_N^{(2)}(\tpi_3 A)}^{a_{2}-1}{Z_N^{(1)}(\tpi_2 A)}^{a_{1}-1}e^{S_N f(x_A)}
\end{multlined} \\
&=
\begin{multlined}[t][15cm]
\frac{1}{Z_N} \sum_{A_r \in \flo{A}^{(r)}_N} {Z_N^{(r-1)}(A_r)}^{a_{r-1}-1} \hspace{-10pt} \sum_{A_{r-1} \in \flo{A}^{(r-1)}_N(A_r)} {Z_N^{(r-2)}(A_{r-1})}^{a_{r-2}-1} \\
\cdots \sum_{A_3 \in \flo{A}^{(3)}_N(A_4)} Z_N^{(2)}(A_3)^{a_2-1} \sum_{A_2 \in \flo{A}^{(2)}_N(A_3)} {Z_N^{(1)}(A_2)}^{a_{1}-1} \underbrace{\sum_{A_1 \in \flo{A}^{(1)}_N(A_2)} e^{S_N f(x_{A_1})}}_{= Z_N^{(1)}(A_2)}
\end{multlined} \\[12pt]
&=
\begin{multlined}[t][11cm]
\frac{1}{Z_N} \sum_{A_r \in \flo{A}^{(r)}_N} {Z_N^{(r-1)}(A_r)}^{a_{r-1}-1} \hspace{-10pt} \sum_{A_{r-1} \in \flo{A}^{(r-1)}_N(A_r)} {Z_N^{(r-2)}(A_{r-1})}^{a_{r-2}-1} \\
\cdots \sum_{A_3 \in \flo{A}^{(3)}_N(A_4)} Z_N^{(2)}(A_3)^{a_2-1} \underbrace{\sum_{A_2 \in \flo{A}^{(2)}_N(A_3)} {Z_N^{(1)}(A_2)}^{a_{1}}}_{= Z_N^{(2)}(A_3)}
\end{multlined} \\
&= \cdots = \frac{1}{Z_N} \sum_{A_r \in \flo{A}^{(r)}_N} {Z_N^{(j-1)}(A_r)}^{a_{r-1}} = 1.
\end{align*}

Although $\sigma_N^{}$ is not generally $T_1$-invariant, the following well-known trick allows us to create a $T_1$-invariant measure $\mu$. We begin by setting
\[ \mu_N^{} = \frac{1}{N} \sum_{k=0}^{N-1} {{T_1}^k}_* \sigma_N^{}. \]
Since $X_1$ is compact, we can take a sub-sequence of $(\mu_N^{})_N$ so that it weakly converges to a probability measure $\mu$ on $X_1$. Then $\mu$ is $T_1$-invariant by the definition of $\mu_N$. We will show that this $\mu$ satisfies 
\[ \sum_{i=1}^r w_i h_{{\pi^{(i-1)}}_* \mu} (T_i, \spa \flo{A}^{(i)}) + w_1 \int_{X_1} f d\mu \geq \lim_{N \to \infty} \frac{\log Z_N}{N}. \]

We first prove
\[ \sum_{i=1}^r w_i H_{{\pi^{(i-1)}}_* \sigma_N} (\flo{A}^{(i)}_N) + w_1 \int_{X_1} \! S_N f d\mu = \log Z_N. \]
To simplify the notations, let 
\begin{align*}
\sigma_N^{(i)} 
&= {\pi^{(i-1)}}_* \sigma_N^{} \\
&= \frac{1}{Z_N} \sum_{B \in \flo{A}^{(1)}_N}{Z_N^{(r-1)}(\tpi_r B)}^{a_{r-1}-1} \cdots {Z_N^{(1)}(\tpi_2 B)}^{a_{1}-1}e^{S_N f(x_B)} \delta_{\pi^{(i)} (x_B)}
\end{align*}
and 
\[ W_N^{(j)} = \sum_{A \in \flo{A}^{(j+1)}_N} {Z_N^{(r-1)}(\tpi_r A)}^{a_{r-1}-1} \cdots {Z_N^{(j+1)}(\tpi_{j+2} A)}^{a_{j+1}-1}{Z_N^{(j)}(A)}^{a_{j}} \log{\left( Z_N^{(j)}(A) \right)}. \]
\begin{claim} \label{claim: h sigma}
We have the following equations:
\begin{align*}
H_{\sigma_N^{}}(\flo{A}^{(1)}_N) &= \spa \log{Z_N} - \int_{X_1} \! S_N f d\sigma_N^{} -\sum_{j=1}^{r-1} \frac{a_{j}-1}{Z_n} W_N^{(j)}, \\[12pt]
H_{\sigma_N^{(i)}}(\flo{A}^{(i)}_N) &= \spa \log{Z_N} - \frac{a_{i-1}}{Z_n} W_N^{(i-1)} - \sum_{j=i}^{r-1} \frac{a_{j}-1}{Z_n} W_N^{(j)} \hspace{8pt} (\hspace{2pt} \it{for} \hspace{4pt} 2 \leq i \leq r \hspace{2pt} ).
\end{align*}
Here, $\sum_{j=r}^{r-1} \frac{a_{j}-1}{Z_n} W_N^{(j)}$ is defined to be $0$.
\end{claim}
\begin{proof}
Let $A \in \flo{A}^{(1)}_N$. We have
\[ \sigma_N^{}(A) 
= \frac{1}{Z_N} {Z_N^{(r-1)}(\tpi_r A)}^{a_{r-1}-1} \cdots {Z_N^{(1)}(\tpi_2 A)}^{a_{1}-1}e^{S_N f(x_A)}. \]
Then
\pagebreak
\begin{flalign*}
& H_{\sigma_N^{}}(\flo{A}^{(1)}_N) = - \sum_{A \in \flo{A}^{(1)}_N} \sigma_N^{}(A) \log{(\sigma_N^{}(A))} &
\end{flalign*} \\[-30pt]
\begin{flalign*}
&
\begin{multlined}[t][11.5cm]
= \spa \log{Z_N} - \underbrace{\sum_{A \in \flo{A}^{(1)}_N} \sigma_N^{}(A) S_N f(x_A)}_{(\mathrm{I})} \\
- \sum_{j=1}^{r-1} \frac{a_j-1}{Z_N} \underbrace{\sum_{A \in \flo{A}^{(1)}_N} {Z_N^{(r-1)}(\tpi_r A)}^{a_{r-1}-1} \cdots {Z_N^{(1)}(\tpi_2 A)}^{a_{1}-1} e^{S_N f(x_A)} \log{\left(Z_N^{(j)}(\tpi_{j+1} A) \right)}}_{(\mathrm{I}\hspace{-0.5pt}\mathrm{I})}.
\end{multlined}
&
\end{flalign*}
For $(\mathrm{I})$, we have 
\begin{align*}
\int_{X_1} \! S_N f d\sigma_N^{}
&= \frac{1}{Z_N} \sum_{A \in \flo{A}^{(1)}_N} {Z_N^{(r-1)}(\tpi_{r}A)}^{a_{r-1}-1} \cdots {Z_N^{(2)}(\tpi_3 A)}^{a_{2}-1}{Z_N^{(1)}(\tpi_2 A)}^{a_{1}-1}e^{S_N f(x_A)} S_N f(x_A) \\
&= (\mathrm{I}).
\end{align*}
We will show that $(\mathrm{I}\hspace{-0.5pt}\mathrm{I}) = W_N^{(j)}$. Let $A' \in \flo{A}^{(j+1)}_N$. Then any $A \in \flo{A}^{(1)}_N(A')$ satisfies $\tpi_{j+1} A = A'$. Hence,
\begin{align*}
(\mathrm{I}\hspace{-0.5pt}\mathrm{I})
&= \sum_{A' \in \flo{A}^{(j+1)}_N} \sum_{A \in \flo{A}^{(1)}_N(A')} {Z_N^{(r-1)}(\tpi_r A)}^{a_{r-1}-1} \cdots {Z_N^{(1)}(\tpi_2 A)}^{a_{1}-1}e^{S_N f(x_A)} \log{\left(Z_N^{(j)}(\tpi_{j+1} A) \right)} \\
&=
\begin{multlined}[t][8cm]
\sum_{A' \in \flo{A}^{(j+1)}_N} {Z_N^{(r-1)}(\tpi_r A')}^{a_{r-1}-1}
\cdots {Z_N^{(j+1)}(\tpi_{j+2} A')}^{a_{j+1}-1} {Z_N^{(j)}(A')}^{a_{j}-1} \log{\left(Z_N^{(j)} (A') \right)} \\
\times \underbrace{\sum_{A \in \flo{A}^{(1)}_N(A')} {Z_N^{(j-1)}(\tpi_j A)}^{a_{j-1}-1} \cdots {Z_N^{(1)}(\tpi_2 A)}^{a_{1}-1}e^{S_N f(x_A)}.}_{(\mathrm{I}\hspace{-0.5pt}\mathrm{I})'}
\end{multlined}
\end{align*}
The term $(\mathrm{I}\hspace{-0.5pt}\mathrm{I})'$ can be calculated similarly to how we showed $\sigma_N^{}(X_1)=1$. Namely,
\begin{align*}
(\mathrm{I}\hspace{-0.5pt}\mathrm{I})'
&=
\begin{multlined}[t][15cm]
\sum_{A_j \in \flo{A}^{(j)}_N(A')} {Z_N^{(j-1)}(A_j)}^{a_{j-1}-1} \hspace{-10pt} \sum_{A_{j-1} \in \flo{A}^{(j-1)}_N(A_j)} {Z_N^{(j-2)}(A_{j-1})}^{a_{j-2}-1} \\
\cdots \sum_{A_3 \in \flo{A}^{(3)}_N(A_4)} Z_N^{(2)}(A_3)^{a_2-1} \sum_{A_2 \in \flo{A}^{(2)}_N(A_3)} {Z_N^{(1)}(A_2)}^{a_{1}-1} \underbrace{\sum_{A_1 \in \flo{A}^{(1)}_N(A_2)} e^{S_N f(x_{A_1})}}_{= Z_N^{(1)}(A_2)}
\end{multlined}
\end{align*}
\begin{flalign*}
& \hspace{20pt} = \cdots = \sum_{A_j \in \flo{A}^{(j)}_N(A')} {Z_N^{(j-1)}(A_j)}^{a_{j-1}} = Z^{(j)}_N(A'). &
\end{flalign*}
Thus, we get
\begin{align*}
(\mathrm{I}\hspace{-0.5pt}\mathrm{I})
&= \sum_{A \in \flo{A}^{(j+1)}_N} {Z_N^{(r-1)}(\tpi_r A)}^{a_{r-1}-1}
\cdots {Z_N^{(j+1)}(\tpi_{j+2} A)}^{a_{j+1}-1}
\cdot {Z_N^{(j)}(A)}^{a_{j}} \log{\left( Z_N^{(j)}(A) \right)} \\[4pt]
&= W_N^{(j)}.
\end{align*}
This completes the proof of the first assertion.

Next, let $2 \leq i \leq r$. For any $A \in \flo{A}^{(i)}_N$, 
\begin{align*}
\sigma_N^{(i)}(A) 
&= \frac{1}{Z_n} \sum_{\substack{B \in \flo{A}^{(1)}_N, \\ \pi^{(i)} (x_B) \in A}} {Z_N^{(r-1)}(\tpi_r B)}^{a_{r-1}-1} \cdots {Z_N^{(1)}(\tpi_2 B)}^{a_{1}-1}e^{S_N f(x_B)} \\
&= 
\begin{multlined}[t][11.5cm]
\frac{1}{Z_n} {Z_N^{(r-1)}(\tpi_r A)}^{a_{r-1}-1} \cdots {Z_N^{(i-1)}(\tpi_i A)}^{a_{i-1}-1} \\
\times \sum_{B \in \flo{A}^{(1)}_N(A)}
{Z_N^{(i-2)}(\tpi_{i-1} B)}^{a_{i-2}-1} \cdots {Z_N^{(1)}(\tpi_2 B)}^{a_{1}-1}e^{S_N f(x_B)}.
\end{multlined}
\end{align*}
As in the evaluation of $(\mathrm{I}\hspace{-0.5pt}\mathrm{I})'$, we have
\[ \sum_{B \in \flo{A}^{(1)}_N(A)} {Z_N^{(i-2)}(\tpi_{i-1} B)}^{a_{i-2}-1} \cdots {Z_N^{(1)}(\tpi_2 B)}^{a_{1}-1}e^{S_N f(x_B)} = {Z^{(i-1)}_N(A)}^{a_{i-1}}. \]
Hence,
\begin{align*}
\sigma_N^{(i)}(A)  &= \frac{1}{Z_n} {Z_N^{(r-1)}(\tpi_r A)}^{a_{r-1}-1} \cdots {Z_N^{(i)}(\tpi_{i+1} A)}^{a_{i}-1} {Z_N^{(i-1)}(A)}^{a_{i-1}}.
\end{align*}
Therefore,
\begin{flalign*}
& \hspace{18pt} H_{\sigma_N^{(i)}}(\flo{A}^{(i)}_N) = - \sum_{A \in \flo{A}^{(i)}_N} \sigma_N^{(i)}(A) \log{\sigma_N^{(i)}(A)} &
\end{flalign*} \\[-35pt]
\begin{align*}
&=
\begin{multlined}[t][15cm]
\log{Z_N} - \frac{1}{Z_n} \sum_{A \in \flo{A}^{(i)}_N} {Z_N^{(r-1)}(\tpi_r A)}^{a_{r-1}-1} \cdots {Z_N^{(i)}(\tpi_{i+1} A)}^{a_{i}-1} {Z_N^{(i-1)}(A)}^{a_{i-1}} \\
\times \log{\left( {Z_N^{(r-1)}(\tpi_r A)}^{a_{r-1}-1} \cdots {Z_N^{(i)}(\tpi_{i+1} A)}^{a_{i}-1} {Z_N^{(i-1)}(A)}^{a_{i-1}} \right)}
\end{multlined} \\[12pt]
&=
\begin{multlined}[t][16cm]
\log{Z_N} - \frac{a_{i-1}}{Z_n} \sum_{A \in \flo{A}^{(i)}_N} {Z_N^{(r-1)}(\tpi_r A)}^{a_{r-1}-1} \cdots {Z_N^{(i)}(\tpi_{i+1} A)}^{a_{i}-1}
{Z_N^{(i-1)}(A)}^{a_{i-1}} \log{\left( Z_N^{(i-1)}(A) \right)} \\
- \sum_{j=i}^{r-1} \frac{a_{j}-1}{Z_n} \sum_{A \in \flo{A}^{(i)}_N} {Z_N^{(r-1)}(\tpi_r A)}^{a_{r-1}-1} \cdots {Z_N^{(i)}(\tpi_{i+1} A)}^{a_{i}-1} {Z_N^{(i-1)}(A)}^{a_{i-1}} \log{\left( Z_N^{(j)}(\tpi_{j+1} A) \right)}.
\end{multlined}
\end{align*}
Note that we have
\pagebreak
\begin{flalign*}
& \sum_{A \in \flo{A}^{(i)}_N} {Z_N^{(r-1)}(\tpi_r A)}^{a_{r-1}-1} \cdots {Z_N^{(i)}(\tpi_{i+1} A)}^{a_{i}-1} {Z_N^{(i-1)}(A)}^{a_{i-1}} \log{\left( Z_N^{(j)}(\tpi_{j+1} A) \right)} & \\[-18pt]
\end{flalign*}
\begin{align*}
&=
\begin{multlined}[t][3cm]
\sum_{A_{j+1} \in \flo{A}^{(j+1)}_N} {Z_N^{(r-1)}(\tpi_r A_{j+1})}^{a_{r-1}-1} \cdots {Z_N^{(j+1)}(\tpi_{j+2} A_{j+1})}^{a_j-1} {Z_N^{(j)}(A_{j+1})}^{a_{j-1}-1} \log{\left( Z_N^{(j)}(A_{j+1}) \right)} \\
\times \sum_{A_j \in \flo{A}^{(j)}_N(A_{j+1})} {Z^{(j-1)}_N(A_j)}^{a_{j-2}-1} \cdots \sum_{A_{i+1} \in \flo{A}^{(i+1)}_N(A_{i+2})} {Z^{(i)}_N(A_{i+1})}^{a_{i+1}-1} \underbrace{\sum_{A_i \in \flo{A}^{(i)}_N(A_{i+1})} {Z^{(i-1)}_N(A_i)}^{a_{i-1}}}_{= Z^{(i)}_N(A_{i+1})}
\end{multlined} \\
&= \cdots = \sum_{A_{j+1} \in \flo{A}^{(j+1)}_N} {Z_N^{(r-1)}(\tpi_r A_{j+1})}^{a_{r-1}-1} \cdots {Z_N^{(j+1)}(\tpi_{j+2} A_{j+1})}^{a_j-1} {Z_N^{(j)}(A_{j+1})}^{a_{j-1}} \log{\left( Z_N^{(j)}(A_{j+1}) \right)}.
\end{align*}
We conclude that
\begin{equation*}
H_{\sigma_N^{(i)}}(\flo{A}^{(i)}_N) = \spa \log{Z_N} - \frac{a_{i-1}}{Z_n} W_N^{(i-1)} - \sum_{j=i}^{r-1} \frac{a_{j}-1}{Z_n} W_N^{(j)}.
\end{equation*}
This completes the proof of the claim.
\end{proof}
By this claim,
\begin{align*}
	\sum_{i=1}^r w_i H_{\sigma_N^{(i)}} (\flo{A}^{(i)}_N) + w_1 \int_{X_1} \! S_N f d\mu
	&= \spa \log{Z_N} - \sum_{i=2}^r \frac{w_i a_{i-1}}{Z_n} W_N^{(i-1)} - \sum_{i=1}^{r-1} \sum_{j=i}^{r-1} \frac{w_i(a_{j}-1)}{Z_n} W_N^{(j)}.
\end{align*}
However, we have
\[ \sum_{i=2}^r w_i a_{i-1} W_N^{(i-1)} + \sum_{i=1}^{r-1} \sum_{j=i}^{r-1} w_i(a_{j}-1) W_N^{(j)} = 0. \]
Indeed, the coefficient of $W_N^{(k)}$ ($1 \leq k \leq r-1$) is
\begin{align*}
w_{k+1}a_k + (a_k - 1) \sum_{i=1}^k w_i &= w_{k+1}a_k + (a_k - 1) a_k a_{k+1} \cdots a_{r-1} \\
&= a_k \{ w_{k+1} - (1-a_k) a_{k+1} a_{k+2} \cdots a_{r-1} \} = 0.
\end{align*}
Thus, we have
\begin{equation} \label{equation: logzn}
\sum_{i=1}^r w_i H_{\sigma_N^{(i)}} (\flo{A}^{(i)}_N) + w_1 \int_{X_1} \! S_N f d\mu = \log{Z_N}.
\end{equation}

Let $\mu^{(i)} = {\pi^{(i-1)}}_* \mu$ and $\mu^{(i)}_N = {\pi^{(i-1)}}_* \mu_N^{}$. 
\begin{lemma} \label{lemma: mu sigma}
Let $N$ and $M$ be natural numbers. For any $1 \leq i \leq r$,
\[ \frac{1}{M} H_{\mu^{(i)}_N} (\flo{A}^{(i)}_M) \geq \frac{1}{N} H_{\sigma_N^{(i)}}(\flo{A}^{(i)}_N) - \frac{2M \log{\abs{\flo{A}^{(i)}}}}{N}. \]
Here, $\abs{\flo{A}^{(i)}}$ is the number of elements in $\flo{A}^{(i)}$.
\end{lemma}
Suppose this is true, and let $N$ and $M$ be natural numbers. Together with \eqref{equation: logzn}, we obtain the following evaluation;
\begin{align*}
\sum_{i=1}^r \frac{w_i}{M} H_{\mu^{(i)}_N}(\flo{A}^{(i)}_M) + w_1 \int_{X_1} \! f d\mu_N^{}
& \geq \sum_{i=1}^r \frac{w_i}{N} H_{\sigma_N^{(i)}}(\flo{A}^{(i)}_N) -\sum_{i=1}^r \frac{2M \log{\abs{\flo{A}^{(i)}}}}{N} + \frac{w_1}{N} \int_{X_1} \! S_N f d\sigma_N^{} \\
&= \frac{\log{Z_N}}{N} - \sum_{i=1}^r \frac{2M \log{\abs{\flo{A}^{(i)}}}}{N}.
\end{align*}
Let $N = N_k \to \infty$ along the sub-sequence $(N_k)$ for which $\mu_{N_k}^{} \rightharpoonup \mu$. This yields
\[ \sum_{i=1}^r \frac{w_i}{M} H_{\mu^{(i)}}(\flo{A}^{(i)}_M) + w_1 \int_{X_1} \! f d\mu \geq \lim_{N \to \infty} \frac{\log{Z_N}}{N}. \]
We let $M \to \infty$ and get
\[ \sum_{i=1}^r w_i h_{\mu^{(i)}}(T_i, \flo{A}^{(i)}) + w_1 \int_{X_1} \! f d\mu \geq \lim_{N \to \infty} \frac{\log{Z_N}}{N}. \]
Hence,
\[ P^{\boldsymbol{a}}_\mathrm{var}(f) \geq P^{\boldsymbol{a}}(f). \]

We are left to prove Lemma \ref{lemma: mu sigma}.
\begin{proof}[Proof of Lemma \ref{lemma: mu sigma}]
This statement appears in the proof of variational principle in \cite[Theorem 8.6]{Walters}, and Tsukamoto also proves it in \cite[Claim 6.3]{Tsukamoto}. The following proof is taken from the latter. We will explain for $i=1$; the same argument works for all $i$.

Let $\flo{A} = \flo{A}^{(1)}$. Recall that $\mu_N^{} = \frac{1}{N} \sum_{k=0}^{N-1} {{T_1}^k}_* \sigma_N^{}$. Since the entropy function is concave (Lemma \ref{lemma: properties of entropy}), we have
\[ H_{\mu_N^{}}(\flo{A}_M) \geq \frac{1}{N} \sum_{k=0}^{N-1} H_{{{T_1}^k}_* \sigma_N^{}}(\flo{A}_M) = \frac{1}{N} \sum_{k=0}^{N-1} H_{\sigma_N^{}}(T_1^{-k}\flo{A}_M). \]
Let $N = qM + r$ with $0 \leq r < M$, then

\begin{align}
\sum_{k=0}^{N-1} H_{\sigma_N^{}}(T_1^{-k}\flo{A}_M)
&= \sum_{s=0}^q \sum_{t=0}^{M-1} H_{\sigma_N^{}}(T_1^{-sM-t}\flo{A}_M) - \sum_{k=N}^{qM+M-1} H_{\sigma_N^{}}(T_1^{-k}\flo{A}_M) \nonumber \\
&\geq \sum_{t=0}^{M-1} \sum_{s=0}^q H_{\sigma_N^{}}(T_1^{-sM-t}\flo{A}_M) - M \log{\abs{\flo{A}_M}} \nonumber \\
&\geq \sum_{t=0}^{M-1} \sum_{s=0}^q H_{\sigma_N^{}}(T_1^{-sM-t}\flo{A}_M) - M^2 \log{\abs{\flo{A}}}. \label{inequality: proof of inequality for entropy of weak limit}
\end{align}
We will evaluate $\sum_{s=0}^q H_{\sigma_N^{}}(T_1^{-sM-t}\flo{A}_M)$ from below for each $0 \leq t \leq M-1$. First, observe that
\[ T_1^{-sM-t}\flo{A}_M = \bigvee_{j=0}^{M-1} T_1^{-sM-t-j}\flo{A}. \]
We have
\[ \{ sM + t + j \spa | \spa 0 \leq s \leq q, 0 \leq j \leq M - 1 \} = \{ t, t + 1, \ldots, t + qM + M - 1 \} \]
without multiplicity. Therefore,
\begin{align*}
H_{\sigma_N^{}}(\flo{A}_N) &\leq H_{\sigma_N^{}}\left( \bigvee_{k=0}^{t+(q+1)M-1} T_1^{-k} \flo{A} \right) \qquad \text{by $N <t+(q+1)M$} \\
&\leq \sum_{s=0}^q H_{\sigma_N^{}}(T_1^{-sM-t}\flo{A}_M) + \sum_{k=0}^{t-1}H_{\sigma_N^{}}(T_1^{-k}\flo{A}) \qquad \text{by Lemma \ref{lemma: properties of entropy}}.
\end{align*}
This implies
\begin{align*}
\sum_{s=0}^q H_{\sigma_N^{}}(T_1^{-sM-t}\flo{A}_M)
&\geq H_{\sigma_N^{}}(\flo{A}_N) - \sum_{k=0}^{t-1}H_{\sigma_N^{}}(T_1^{-k}\flo{A}) \\
&\geq H_{\sigma_N^{}}(\flo{A}_N) - M \log{\abs{\flo{A}}} \qquad \text{by $t < M$}.
\end{align*}
Now, we sum over $t$ and obtain
\begin{align*}
\sum_{t=1}^{M-1} \sum_{s=0}^q H_{\sigma_N^{}}(T_1^{-sM-t}\flo{A}_M)
&\geq M H_{\sigma_N^{}}(\flo{A}_N) - M^2 \log{\abs{\flo{A}}}.
\end{align*}
Combining with \eqref{inequality: proof of inequality for entropy of weak limit}, this implies
\begin{equation*}
\sum_{k=0}^{N-1} H_{\sigma_N^{}}(T_1^{-k}\flo{A}_M)
\geq M H_{\sigma_N^{}}(\flo{A}_N) - 2M^2 \log{\abs{\flo{A}}}.
\end{equation*}
It follows that
\begin{equation*}
\frac{1}{M} H_{\mu_N^{}}(\flo{A}_N) \geq \frac{1}{MN} \sum_{k=0}^{N-1} H_{\sigma_N^{}}(T_1^{-k}\flo{A}_M) \geq \frac{1}{N} H_{\sigma_N^{}}(\flo{A}_N) - \frac{2M \log{\abs{\flo{A}}}}{N}.
\end{equation*}
\end{proof}
This completes the proof of Theorem \ref{theorem: first half of the main theorem}.
\end{proof}

\section{\texorpdfstring{Proof of $P^{\boldsymbol{a}}_\mathrm{var}(f) \leq P^{\boldsymbol{a}}(f)$.}%
                               {subtitle}} \label{section: proof of Pvar is smaller}
It seems difficult to implement the zero-dimensional trick to prove $P^{\boldsymbol{a}}_\mathrm{var}(f) \leq P^{\boldsymbol{a}}(f)$. Hence, the proof is more complicated.
\begin{theorem} Suppose that $(X_i, \spa T_i)$ ($i=1, \spa 2, \spa \ldots, \spa r$) are dynamical systems and $\pi_i: X_i \rightarrow X_{i+1} \hspace{5pt} (i=1, \spa 2, \spa ... \spa , \spa r-1)$ are factor maps. Then we have
\[ P^{\boldsymbol{a}}_\mathrm{var}(f) \leq P^{\boldsymbol{a}}(f) \]
for any continuous function $f: X_1 \rightarrow \mathbb{R}$.
\end{theorem}
\begin{proof}
Take and fix $\mu \in \flo{M}^{T_1}(X_1)$. Let $\mu_i= {\pi^{(i-1)}}_* \mu$. We need to prove
\[ \sum_{i=1}^r w_i h_{\mu_i}(T_i) + w_1 \int_{X_1} \! f d\mu \leq P^{\boldsymbol{a}}(f, \boldsymbol{T}). \]
However, the following argument assures that giving an evaluation up to a constant is sufficient: suppose there is a positive number $C$ which does not depend on $f$ nor $(T_i)_i$ satisfying
\begin{equation} \label{equation: up to constant}
\sum_{i=1}^r w_i h_{\mu_i}(T_i) + w_1 \int_{X_1} \! f d\mu \leq P^{\boldsymbol{a}}(f, \boldsymbol{T}) + C.
\end{equation}
Applying this to $S_mf$ and $\boldsymbol{T}^m = ({T_i}^m)_i$ for $m \in \mathbb{N}$ yields
\[ \sum_{i=1}^r w_i h_{\mu_i}({T_i}^m) + w_1 \int_{X_1} \! S_mf d\mu \leq P^{\boldsymbol{a}}(S_mf, \boldsymbol{T}^m) + C. \]
We employ Lemma \ref{lemma: multiplication penetrates} and get
\[ m \sum_{i=1}^r w_i h_{\mu_i}(T_i) + m w_1 \int_{X_1} \! f d\mu \leq mP^{\boldsymbol{a}}(f, \boldsymbol{T}) + C. \]
Divide by $m$ and let $m \to \infty$. We obtain the desired inequality
\[ \sum_{i=1}^r w_i h_{\mu_i}(T_i) + w_1 \int_{X_1} \! f d\mu \leq P^{\boldsymbol{a}}(f, \boldsymbol{T}). \]
Therefore, we only need to prove (\refeq{equation: up to constant}).

Let $\flo{A}^{(i)} = \{ A^{(i)}_1, A^{(i)}_2, \cdots, A^{(i)}_{m_i} \}$ be an arbitrary partition of $X_i$ for each $i$. We will prove
\begin{equation*}
\sum_{i=1}^r w_i h_{\mu_i}(T_i, \flo{A}^{(i)}) + w_1 \int_{X_1} \! f d\mu \leq P^{\boldsymbol{a}}(f, \boldsymbol{T}) + C.
\end{equation*}
We start by approximating elements of $\flo{A}^{(i)}$ with compact sets using backward induction. For $1\leq i \leq r$, let
\[ \Lambda_i^{0} = \{0, 1, \cdots, m_r\} \times \{0, 1, \cdots, m_{r-1}\} \times \cdots \times \{0, 1, \cdots, m_{i+1}\} \times \{0, 1, \cdots, m_i\},\]
\[ \Lambda_i = \{0, 1, \cdots, m_r\} \times \{0, 1, \cdots, m_{r-1}\} \times \cdots \times \{0, 1, \cdots, m_{i+1}\} \times \{1, 2, \cdots, m_i\}. \]
We will denote an element $(j_r, j_{r-1}, \cdots, j_i)$ in $\Lambda_i^{0}$ or $\Lambda_i$ by $j_r j_{r-1} \cdots j_i$.
For each $A^{(r)}_j \in \flo{A}^{(r)}$, take a compact set $C^{(r)}_j \subset A^{(r)}_j$ such that 
\[ \log{m_r} \cdot \sum_{j=1}^{m_r} \mu_r(A^{(r)}_j \setminus C^{(r)}_j) < 1. \]
Define $C^{(r)}_0$ as the remainder of $X_r$, which may not be compact;
\[ C^{(r)}_0 = \bigcup_{j=1}^{m_r} A^{(r)}_j \setminus C^{(r)}_j = X_r \setminus \bigcup_{j=1}^{m_r} C^{(r)}_j. \]
Then $\flo{C}^{(r)} := \{ C^{(r)}_0, C^{(r)}_1, \cdots, C^{(r)}_{m_r} \}$ is a measurable partition of $X_r$.

Next, consider the partition $\pi_{r-1}^{-1}(\flo{C}^{(r)}) \vee \flo{A}^{(r-1)}$ of $X_{r-1}$. For $j_r j_{r-1} \in \Lambda_{r-1}$, let
\[ B^{(r-1)}_{j_r j_{r-1}} = \pi_{r-1}^{-1}(C^{(r)}_{j_r}) \cap A^{(r-1)}_{j_{r-1}}. \]
Then
\[ \pi_{r-1}^{-1}(\flo{C}^{(r)}) \vee \flo{A}^{(r-1)} = \left\{ B^{(r-1)}_{j_r j_{r-1}}  \spa \middle| \spa \text{ $j_r j_{r-1} \in \Lambda_{r-1}$ } \right\}, \]
and for each $j_r \in \Lambda_r^0$
\[ \bigcup_{j_{r-1}=1}^{m_{r-1}} B^{(r-1)}_{j_r j_{r-1}} = \pi_{r-1}^{-1}(C^{(r-1)}_{j_r}). \]
For each $j_rj_{r-1} \in \Lambda_{r-1}$, take a compact set $C^{(r-1)}_{j_r j_{r-1}} \subset B^{(r-1)}_{j_r j_{r-1}}$ (which could be empty) such that 
\[ \log{\abs{\Lambda_{r-1}}} \cdot \sum_{j_r = 0}^{m_r} \sum_{j_{r-1}=1}^{m_{r-1}} \mu_{r-1}(B^{(r-1)}_{j_r j_{r-1}} \setminus C^{(r-1)}_{j_r j_{r-1}}) < 1. \]
Define $C^{(r-1)}_{j_r 0}$ as the remainder of $\pi_{r-1}^{-1}(C^{(r)}_{j_r})$;
\[ C^{(r-1)}_{j_r 0} = \pi_{r-1}^{-1}(C^{(r)}_{j_r}) \setminus \bigcup_{j_{r-1}=1}^{m_{r-1}} C^{(r-1)}_{j_r j_{r-1}}. \]
Then $\flo{C}^{(r-1)} = \left\{ C^{(r-1)}_{j_r j_{r-1}} \spa \middle| \spa j_r j_{r-1} \in \Lambda_{r-1}^{0} \right\} $ is a measurable partition of $X_{r-1}$.

Continue in this manner, and suppose we have obtained the partition $\flo{C}^{(k)} =  \left\{ C^{(k)}_J \spa \middle| \spa J \in \Lambda_k^{0} \right\}$ of $X_k$ for $k = i+1, i+2, \ldots, r$.
We will define $\flo{C}^{(i)}$. Each element in $\pi_i^{-1}(\flo{C}^{(i+1)}) \vee \flo{A}^{(i)}$ can be expressed using $J' \in \Lambda_{i+1}^0$ and $j_i \in \{1,2, \ldots, m_i\}$ by
\[ B_{J'j_i}^{(i)} = \pi_i^{-1}(C_{J'}^{(i+1)}) \cap A_{j_i}^{(i)}. \]
Choose a compact set $C_J^{(i)} \subset B_J^{(i)}$ for each $J \in \Lambda_{i}$ so that
\[ \log{\abs{\Lambda_i}} \cdot \sum_{J' \in \Lambda_{i+1}^0} \sum_{j_i=1}^{m_i} \mu_i \left( B^{(i)}_{J'j_i} \setminus C^{(i)}_{J'j_i} \right) < 1. \]
Finally, for $J' \in \Lambda_{j+1}^{0}$, let
\[ C_{J'0}^{(i)} = \pi_i^{-1}(C_{J'}^{(i+1)}) \setminus \bigcup_{j_i=1}^{m_i} C_{J'j_i}^{(i)}. \]
Set $\flo{C}^{(i)} = \left\{ C^{(i)}_J \spa \middle| \spa J \in \Lambda_i^{0} \right\}$.
This is a partition of $X_i$.
\begin{lemma} \label{lemma: compact approximation}
For $\flo{C}^{(i)}$ constructed above, we have
\[ h_{\mu_i}(T_i, \flo{A}^{(i)}) \leq h_{\mu_i}(T_i, \flo{C}^{(i)}) + 1. \]
\end{lemma}
\begin{proof}
By Lemma \ref{lemma: properties of entropy},
\begin{align*}
h_{\mu_i}(T_i, \flo{A}^{(i)})
&\leq h_{\mu_i} \! \left(T_i, \flo{A}^{(i)} \vee \pi_i^{-1}(\flo{C}^{(i+1)})\right) \\
&\leq h_{\mu_i}(T_i, \flo{C}^{(i)}) + H_{\mu_i}\left( \flo{A}^{(i)} \vee \pi_i^{-1}(\flo{C}^{(i+1)}) \spa \middle| \spa \flo{C}^{(i)} \right).
\end{align*}
Since $C^{(i)}_J \subset B^{(i)}_J$ for $J \in \Lambda_i$,
\begin{flalign*}
& H_{\mu_i}\left( \flo{A}^{(i)} \vee \pi_i^{-1}(\flo{C}^{(i+1)}) \spa \middle| \spa \flo{C}^{(i)} \right) &
\end{flalign*}
\begin{align*}
&= - \sum_{\substack{J \in \Lambda_i^0 \\ \mu_i(C_J^{(i)}) \ne 0}} \mu_i(C_J^{(i)}) \sum_{K \in \Lambda_i} \frac{\mu_i(B_{K}^{(i)} \cap C_J^{(i)})}{\mu_i(C_J^{(i)})} \log{\left( \frac{\mu_i(B_{K}^{(i)} \cap C_J^{(i)})}{\mu_i(C_J^{(i)})} \right)} \\
&= - \sum_{\substack{J' \in \Lambda_{i+1}^0 \\ \mu_i(C_{J'0}^{(i)}) \ne 0}} \mu_i(C_{J'0}^{(i)}) \sum_{j_i=1}^{m_i} \frac{\mu_i(B_{J' j_i}^{(i)} \cap C_{J'0}^{(i)})}{\mu_i(C_{J'0}^{(i)})} \log{\left( \frac{\mu_i(B_{J' j_i}^{(i)} \cap C_{J'0}^{(i)})}{\mu_i(C_{J'0}^{(i)})} \right)}.
\end{align*}
By Lemma \ref{lemma: calculus}, we have
\[ - \hspace{3pt} \sum_{j_i=1}^{m_i} \frac{\mu_i(B_{J' j_i}^{(i)} \cap C_{J'0}^{(i)})}{\mu_i(C_{J'0}^{(i)})} \log{\left( \frac{\mu_i(B_{J' j_i}^{(i)} \cap C_{J'0}^{(i)})}{\mu_i(C_{J'0}^{(i)})} \right)} \leq \log{\abs{\Lambda_i}}. \]
Therefore,
\begin{align*}
H_{\mu_i}\left( \flo{A}^{(i)} \vee \pi_i^{-1}(\flo{C}^{(i+1)}) \spa \middle| \spa \flo{C}^{(i)} \right)
&\leq \log{\abs{\Lambda_i}} \sum_{J' \in \Lambda_{i+1}^0} \mu_i\left(\pi_i^{-1}(C_{J'}^{(i+1)}) \setminus \bigcup_{j_i=1}^{m_i} C_{J'j_i}^{(i)}\right) < 1.
\end{align*}
\end{proof}
Recall the definition of $\boldsymbol{w}$ in (\refeq{definition: w}). We have
\begin{flalign*}
& \sum_{i=1}^r w_ih_{\mu_i}(T_i, \flo{C}^{(i)}) +  w_1 \int_{X_1} f d\mu &
\end{flalign*} \\[-35pt]
\begin{align*}
&=
\begin{multlined}[t][15cm]
\lim_{N \to \infty} \frac{1}{N} \Bigg\{ H_{\mu_r}(\flo{C}^{(r)}_N) + a_1a_2 \cdots a_{r-1} N \int_{X_1} f d\mu \\
+ \sum_{i=1}^{r-1} a_i a_{i+1} \cdots a_{r-1}\left(H_{\mu_i}(\flo{C}^{(i)}_N) - H_{\mu_{i+1}}(\flo{C}^{(i+1)}_N)\right) \Bigg\}
\end{multlined} \\
&=
\begin{multlined}[t][14cm]
\lim_{N \to \infty} \frac{1}{N} \Bigg\{ H_{\mu_r}(\flo{C}^{(r)}_N) + a_1a_2 \cdots a_{r-1} \int_{X_1} S_Nf d\mu \\
+ \sum_{i=1}^{r-1} a_i a_{i+1} \cdots a_{r-1} H_{\mu_i}\left(\flo{C}^{(i)}_N \middle| \pi_i^{-1}(\flo{C}^{(i+1)}_N) \right) \Bigg\}.
\end{multlined}
\end{align*}
Here, we used the relation
\begin{align*}
H_{\mu_i}(\flo{C}^{(i)}_N) - H_{\mu_{i+1}}(\flo{C}^{(i+1)}_N)
&= H_{\mu_i}(\flo{C}^{(i)}_N) - H_{\mu_i}(\pi_i^{-1}(\flo{C}^{(i+1)}_N)) \\
&=H_{\mu_i}\left(\flo{C}^{(i)}_N \middle| \pi_i^{-1}(\flo{C}^{(i+1)}_N) \right).
\end{align*}

We fix $N$ and evaluate from above the following terms using backward induction:
\begin{align}
H_{\mu_r}(\flo{C}^{(r)}_N) + a_1a_2 \cdots a_{r-1} \int_{X_1} S_Nf d\mu + \sum_{i=1}^{r-1} a_i a_{i+1} \cdots a_{r-1} H_{\mu_i}\left(\flo{C}^{(i)}_N \middle| \pi_i^{-1}(\flo{C}^{(i+1)}_N) \right). \label{eq: Hvar}
\end{align}
First, consider the term
\[a_1a_2 \cdots a_{r-1} \left( H_{\mu}\left(\flo{C}^{(1)}_N \middle| \pi_1^{-1}(\flo{C}^{(2)}_N) \right) + \int_{X_1} S_Nf d\mu \right). \]
For $C \in \flo{C}^{(i+1)}_N$, let $\flo{C}^{(i)}_N(C) = \{ D \in \flo{C}^{(i)}_N \spa | \spa \pi_i(D) \subset C \}$, then by Lemma \ref{lemma: calculus},
\begin{flalign*}
& H_{\mu}\left(\flo{C}^{(1)}_N \middle| \pi_1^{-1}(\flo{C}^{(2)}_N) \right) + \int_{X_1} S_Nf d\mu & \\
& \hspace{80pt} \leq \sum_{\substack{C \in \flo{C}^{(2)}_N \\ \mu_2(C) \ne 0}} \mu_2(C) \left\{ \sum_{D \in \flo{C}^{(1)}_N(C)} \left( -\frac{\mu(D)}{\mu_2(C)}\log{\frac{\mu(D)}{\mu_2(C)}} + \frac{\mu(D)}{\mu_2(C)} \sup_{D} S_Nf \right) \right\} & \\
& \hspace{80pt} \leq \sum_{C \in \flo{C}^{(2)}_N} \mu_2(C) \log{ \sum_{D \in \flo{C}^{(1)}_N(C)} e^{\sup_DS_Nf}}. &
\end{flalign*}
Applying this inequality to (\refeq{eq: Hvar}), the following term appears:
\begin{equation}
 a_2a_3 \cdots a_{r-1} \left( H_{\mu_2}\left(\flo{C}^{(2)}_N \middle| \pi_2^{-1}(\flo{C}^{(3)}_N) \right) + a_1 \sum_{C \in \flo{C}^{(2)}_N} \mu_2(C) \log{ \sum_{D \in \flo{C}^{(1)}_N(C)} e^{\sup_DS_Nf}} \right).
\end{equation}
This can be evaluated similarly using Lemma \ref{lemma: calculus} as \\
\begin{flalign*}
& H_{\mu_2}\left(\flo{C}^{(2)}_N \middle| \pi_2^{-1}(\flo{C}^{(3)}_N) \right) + a_1 \sum_{C \in \flo{C}^{(2)}_N} \mu_2(C) \log{ \sum_{D \in \flo{C}^{(1)}_N(C)} e^{\sup_DS_Nf}} &
\end{flalign*} \\[-27pt]
\begin{align*}
&= \sum_{\substack{C \in \flo{C}^{(3)}_N \\ \mu_3(C) \ne 0}} \mu_3(C) \left\{ \sum_{D \in \flo{C}^{(2)}_N(C)} \left( -\frac{\mu_2(D)}{\mu_3(C)}\log{\frac{\mu_2(D)}{\mu_3(C)}} + \frac{\mu_2(D)}{\mu_3(C)} \log{ \left( \sum_{E \in \flo{C}^{(1)}_N(D)} e^{\sup_ES_Nf} \right)^{a_1} } \right) \right\} \\
&\leq \sum_{C \in \flo{C}^{(3)}_N} \mu_3(C) \log{ \sum_{D \in \flo{C}^{(2)}_N(C)}\left( \sum_{E \in \flo{C}^{(1)}_N(D)} e^{\sup_ES_Nf} \right)^{a_1}}.
\end{align*}
Continue likewise and obtain the following upper bound for (\refeq{eq: Hvar}):
\begin{equation} \label{eq: ValueVar}
\log{ \sum_{C^{(r)} \in \flo{C}^{(r)}_N} \left( \sum_{C^{(r-1)} \in \flo{C}^{(r-1)}_N(C^{(r)})} \left( \cdots \left( \sum_{C^{(1)} \in \flo{C}^{(1)}_N(C^{(2)})} e^{\sup_{C^{(1)}}S_Nf} \right)^{a_1} \cdots \right)^{a_{r-2}}\right)^{a_{r-1}}}.
\end{equation}

For $1\leq i \leq r$, let $\flo{C}^{(i)}_c = \{ C \in \flo{C}^{(i)} \spa | \spa C \text{ is compact} \}$. There is a positive number $\vep_i$ such that $d^{(i)}(y_1, y_2) > \vep_i$ for any $C_1, C_2 \in \flo{C}^{(i)}_c$ and $y_1 \in C_1, y_2 \in C_2$. Fix a positive number $\vep$ with
\begin{equation} \label{eq: epsilon}
\vep < \min_{1\leq i \leq r} \vep_i.
\end{equation}
Let $\flo{F}^{(i)}$ be a chain of open ($N$, $\vep$)-covers of $X_i$ (see Definition \ref{definition: chain of covers}). Consider
\begin{flalign*}
& \log{ \flo{P}^{\boldsymbol{a}}\left( f, \spa N, \spa \vep, \spa (\flo{F}^{(i)})_i \right) } &
\end{flalign*} \\[-35pt]
\begin{align}
&= \log{ \sum_{U^{(r)} \in \flo{F}^{(r)}} \left( \sum_{U^{(r-1)} \in \flo{F}^{(r-1)}(U^{(r)})} \left( \cdots \left( \sum_{U^{(1)} \in \flo{F}^{(1)}(U^{(2)})} e^{\sup_{U^{(1)}}S_Nf} \right)^{a_1} \cdots \right)^{a_{r-2}}\right)^{a_{r-1}}}. \label{eq: ValuePress}
\end{align}
We will evaluate \eqref{eq: ValueVar} from above by \eqref{eq: ValuePress} up to a constant. We need the next lemma.
\begin{lemma} \label{lemma: combinatorics}
	\begin{enumerate}
	\item[(1)] For any $V \subset X_r$ with $\diam(V, d^{(r)}_N) < \vep$,
\[ \left| \left\{ D \in \flo{C}^{(r)}_N \spa \middle| \spa D \cap V \ne \varnothing \right\} \right| \leq 2^N. \]
	\item[(2)] Let $1\leq i \leq r-1$ and $C \in \flo{C}^{(i+1)}_N$. For any $V \subset X_i$ with $\diam(V, d^{(i)}_N) < \vep$,
\[ \left| \left\{ D \in \flo{C}^{(i)}_N(C) \spa \middle| \spa D \cap V \ne \varnothing \right\} \right| \leq 2^N. \]
	\end{enumerate}
\end{lemma}
\begin{proof}
(1) $D \in \flo{C}^{(r)}_N$ can be expressed using $C^{(r)}_{k_s} \in \flo{C}^{(r)}$ ($s =0, 1, \ldots, N-1$) as
\[ D = C^{(r)}_{k_0^{}} \cap T_r^{-1} C^{(r)}_{k_1^{}} \cap T_r^{-2} C^{(r)}_{k_2^{}} \cap \cdots \cap T_r^{-N+1} C^{(r)}_{k_{N-1}^{}}. \]
If $D \cap V \ne \varnothing$, we have $T_r^{-s}(C^{(r)}_{k_s}) \cap V \ne \varnothing$ for every $0 \leq s \leq N-1$. Then for each $s$
\[ \varnothing \ne T_r^s\left( T_r^{-s}(C^{(r)}_{k_s}) \cap V \right) \subset C^{(r)}_{k_s} \cap T_r^s(V). \]
By \eqref{eq: epsilon}, each $k_s$ is either $0$ or one of the elements in $\{ 1, 2, \ldots, m_r \}$. Therefore, there are at most $2^N$ such sets. \\
(2) The proof works in the same way as in (1). $C$ can be written using $J_k \in \Lambda_{i+1}^0$ ($k=0, 1, \ldots, N-1$) as
\[ C = C^{(i+1)}_{J_0} \cap T_{i+1}^{-1} C^{(i+1)}_{J_1} \cap T_{i+1}^{-2} C^{(i+1)}_{J_2} \cap \cdots \cap T_{i+1}^{-N+1} C^{(i+1)}_{J_{N-1}}. \]
Then any $D \in \flo{C}^{(i)}_N(C)$ is of the form
\[ D = C^{(i)}_{J_0k_0} \cap T_i^{-1}C^{(i)}_{J_1k_1} \cap T_i^{-2}C^{(i)}_{J_2k_2} \cap \cdots \cap T_i^{-N+1}C^{(i)}_{J_{N-1}k_{N-1}} \]
with $0 \leq k_l \leq m_i$ ($l = 1, 2, \ldots, N-1$). If $D \cap V \ne \varnothing$, then each $k_l$ is either $0$ or one of the elements in $\{ 1, 2, \ldots, m_i \}$. Therefore, there are at most $2^N$ such sets.
\end{proof}
For any $C^{(1)} \in \flo{C}^{(1)}_N$, there is $V \in \flo{F}^{(1)}$ with $V \cap C^{(1)} \ne \varnothing$ and
\[ e^{\sup_{C^{(1)}}S_ Nf} \leq e^{\sup_V S_N f}. \]
Let $C^{(2)} \in \flo{C}^{(2)}_N$, then by Lemma \ref{lemma: combinatorics},
\[ \sum_{C^{(1)} \in \flo{C}^{(1)}_N(C^{(2)})} e^{\sup_{C^{(1)}}S_ Nf} \leq \sum_{\substack{U \in \flo{F}^{(2)} \\ U \cap C^{(2)} \ne \varnothing}} 2^N \sum_{V \in \flo{F}^{(1)}(U)} e^{\sup_V S_N f}. \]
By Lemma \ref{lemma: calculus},
\[ \left( \sum_{C^{(1)} \in \flo{C}^{(1)}_N(C^{(2)})} e^{\sup_{C^{(1)}}S_ Nf} \right)^{a_1} \leq 2^{a_1N} \sum_{\substack{U \in \flo{F}^{(2)} \\ U \cap C^{(2)} \ne \varnothing}} \left( \sum_{V \in \flo{F}^{(1)}(U)} e^{\sup_V S_N f} \right)^{a_1}. \]
For $C^{(3)} \in \flo{C}^{(3)}_N$, we apply Lemma \ref{lemma: combinatorics} and Lemma \ref{lemma: calculus} similarly and obtain
\begin{align*}
\begin{multlined}[t][15cm]
\left( \sum_{C^{(2)} \in \flo{C}^{(2)}_N(C^{(3)})} \left( \sum_{C^{(1)} \in \flo{C}^{(1)}_N(C^{(2)})} e^{\sup_{C^{(1)}}S_ Nf} \right)^{a_1} \right)^{a_2} \\
\leq 2^{a_1a_2N} 2^{a_2N} \sum_{\substack{O \in \flo{F}^{(3)} \\ O \cap C^{(3)} \ne \varnothing}} \left( \sum_{U \in \flo{F}^{(2)}(O)} \left( \sum_{V \in \flo{F}^{(1)}(U)} e^{\sup_V S_N f} \right)^{a_1} \right)^{a_2}.
\end{multlined}
\end{align*}
We continue this reasoning and get
\begin{align*}
\begin{multlined}[t][15cm]
\sum_{C^{(r)} \in \flo{C}^{(r)}_N} \left( \sum_{C^{(r-1)} \in \flo{C}^{(r-1)}_N(C^{(r)})} \left( \cdots \left( \sum_{C^{(1)} \in \flo{C}^{(1)}_N(C^{(2)})} e^{\sup_{C^{(1)}}S_Nf} \right)^{a_1} \cdots \right)^{a_{r-2}}\right)^{a_{r-1}} \\
\leq 2^{\alpha N} \sum_{U^{(r)} \in \flo{F}^{(r)}} \left( \sum_{U^{(r-1)} \in \flo{F}^{(r-1)}(U^{(r)})} \left( \cdots \left( \sum_{U^{(1)} \in \flo{F}^{(1)}(U^{(2)})} e^{\sup_{U^{(1)}}S_Nf} \right)^{a_1} \cdots \right)^{a_{r-2}}\right)^{a_{r-1}}.
\end{multlined}
\end{align*}
Here $\alpha$ stands for $\sum_{i=1}^{r-1} a_i a_{i+1}\cdots a_{r-1}$. We take the logarithm of both sides; the left-hand side equals \eqref{eq: ValueVar}, which is an upper bound for \eqref{eq: Hvar}. Furthermore, consider the infimum over the chain of open ($N$, $\vep$)-covers $(\flo{F}^{(i)})_i$. By Remark \ref{remark: chains of covers}, this yields
\begin{align*}
\begin{multlined}[t][23cm]
H_{\mu_r}(\flo{C}^{(r)}_N) + a_1a_2 \cdots a_{r-1} \int_{X_1} S_Nf d\mu + \sum_{i=1}^{r-1} a_i a_{i+1} \cdots a_{r-1} H_{\mu_i}\left(\flo{C}^{(i)}_N \middle| \pi_i^{-1}(\flo{C}^{(i+1)}_N) \right) \\
\leq \log{P^{\boldsymbol{a}}_r(X_r, \spa f, \spa N, \spa \vep)} + \alpha N \log{2}. \\
\end{multlined}
\end{align*}
Divide by $N$, then let $N \to \infty$ and $\vep \to 0$. We obtain
\[ \sum_{i=1}^r w_i h_{\mu_i}(T_i, \flo{C}^{(i)}) +  w_1 \int_{X_1} f d\mu \leq P^{\boldsymbol{a}}(f, \boldsymbol{T}) + \alpha \log{2}. \]
Lemma \ref{lemma: compact approximation} yields
\[ \sum_{i=1}^r w_i h_{\mu_i}(T_i, \flo{A}^{(i)}) +  w_1 \int_{X_1} f d\mu \leq P^{\boldsymbol{a}}(f, \boldsymbol{T}) + \alpha \log{2} + r. \]
We take the supremum over the partitions $(\flo{A}^{(i)})_i$:
\[ \sum_{i=1}^r w_i h_{\mu_i}(T_i) + w_1 \int_{X_1} \! f d\mu \leq P^{\boldsymbol{a}}(f, \boldsymbol{T}) + \alpha \log{2} + r. \]
By the argument at the beginning of this proof, we conclude that
\[ \sum_{i=1}^r w_i h_{\mu_i}(T_i) + w_1 \int_{X_1} \! f d\mu \leq P^{\boldsymbol{a}}(f, \boldsymbol{T}). \]
\end{proof}

\section{Example: Sofic Sets} \label{section: example: sofic set}
Kenyon--Peres \cite{Kenyon--Peres: sofic} calculated the Hausdorff dimension of sofic sets in $\mathbb{T}^2$. In this section, we will see that we can calculate the Hausdorff dimension of certain sofic sets in $\mathbb{T}^d$ with arbitrary $d$. We give an example for the case $d=3$.

\subsection{Definition of Sofic Sets} \label{subsection: definition of sofic sets}
This subsection referred to \cite{Kenyon--Peres: sofic}. Weiss \cite{Weiss} defined {\it{sofic systems}} as subshifts which are factors of shifts of finite type. Boyle, Kitchens, and Marcus proved in \cite{Boyle--Kitchens--Marcus} that this is equivalent to the following definition.
\begin{definition}[{{\cite[Proposition 3.6]{Kenyon--Peres: sofic}}}]
Consider a finite directed graph $G = \langle V, E \rangle$ in which loops and multiple edges are allowed. Suppose its edges are colored in $l$ colors in a ``right-resolving'' fashion: every two edges emanating from the same vertex have different colors. Then the set of color sequences that arise from infinite paths in $G$ is called the \textbf{sofic system}.
\end{definition}

Let $m_1 \leq m_2 \leq \cdots \leq m_r$ be natural numbers, $T$ an endomorphism on $\mathbb{T}^r = \mathbb{R}^r/\mathbb{Z}^r$ represented by the diagonal matrix $A = \mathrm{diag}(m_1, m_2, \ldots, m_r)$, and $D = \prod_{i=1}^r \{0, 1, \ldots, m_i-1\}$. Define a map $R_r: D^{\mathbb{N}} \rightarrow \mathbb{T}^r$ by
\[ R_r((e^{(n)})_{n=1}^{\infty}) = \left( \sum_{k=0}^{\infty} \frac{e^{(k)}_1}{{m_1}^k}, \cdots, \sum_{k=0}^{\infty} \frac{e^{(k)}_r}{{m_r}^k} \right) \]
where $e^{(k)} = (e^{(k)}_1, \cdots, e^{(k)}_r) \in D$ for each $k$. Suppose the edges in some finite directed graph are labeled by the elements in $D$ in the right-resolving fashion, and let $S \subset D^{\mathbb{N}}$ be the resulting sofic system. The image of $S$ under $R_r$ is called a \textbf{sofic set}.

\subsection{An example of a sofic set}

\begin{figure}[h!] \label{figure: a sofic system}
\centering
\makebox[0pt]{\includegraphics[width=15cm]{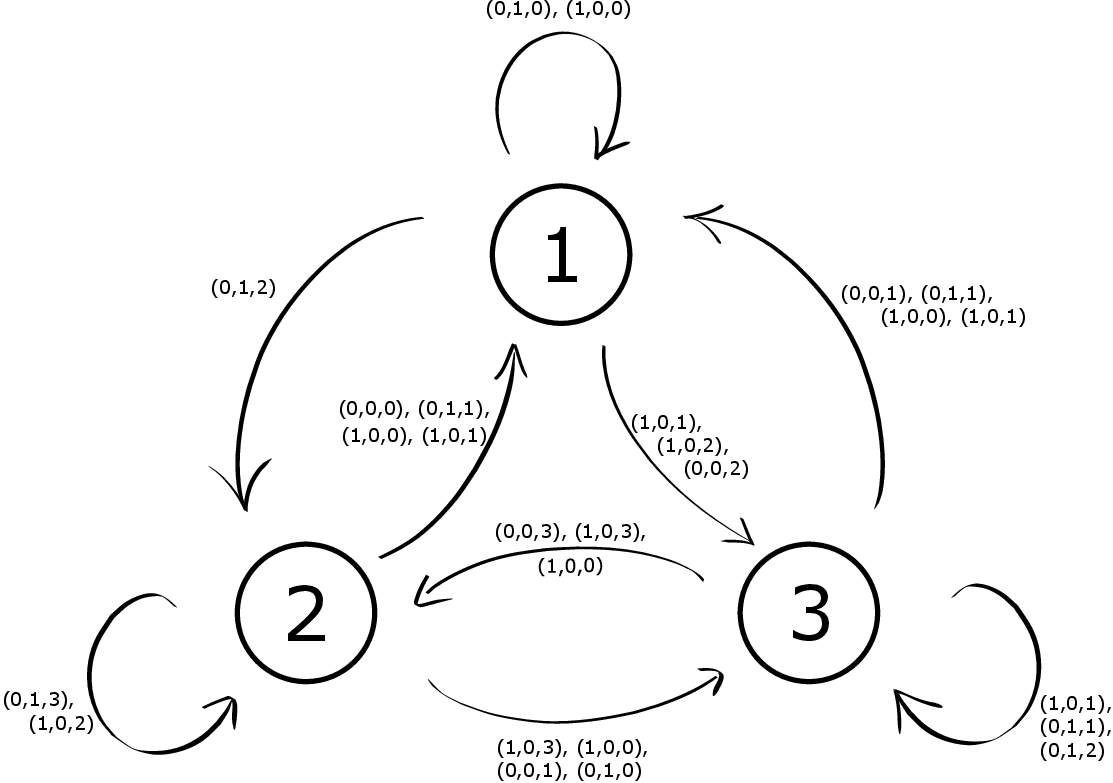}}
\vspace{20pt}%
\caption{Directed graph $G$}
\end{figure}

Here we will look at an example of a sofic set and calculate its Hausdorff dimension via its weighted topological entropy. Let $D = \{0, 1\} \times \{0, 1, 2\} \times \{0, 1, 2, 3\}$ and consider the directed graph $G = \langle V, E \rangle$ with $V = \{1, 2, 3\}$ and $D$-labeled edges in Figure \ref{figure: a sofic system}.

Let $Y_1 \subset D^{\mathbb{N}}$ be the resulting sofic system. Let $C = \{0, 1\} \times \{0, 1, 2\}$ and $B = \{0, 1\}$. Define $p_1: D \rightarrow C$ and $p_2: C \rightarrow B$ by
\[ p_1(i, j, k) = (i, j), \quad p_2(i, j) = i. \]
Let $p_1^{\mathbb{N}}: D^{\mathbb{N}} \rightarrow C^{\mathbb{N}}$ and $p_2^{\mathbb{N}}: C^{\mathbb{N}} \rightarrow B^{\mathbb{N}}$ be the product map of $p_1$ and $p_2$, respectively. Set $Y_2 = p_1^{\mathbb{N}}(Y_1)$ and $Y_3 = p_2^{\mathbb{N}}(Y_2)$. Note that $Y_2 = \{ (0, 0), (1,0), (0, 1) \}^{\mathbb{N}}$ and $Y_3 = \{0, 1\}^{\mathbb{N}}$, meaning they are full shifts.

The sets $X_i = R_i(Y_i)$ $(i = 1, 2, 3)$ are sofic sets. Define $\pi_1: X_1 \rightarrow X_2$ and $\pi_2: X_2 \rightarrow X_3$ by
\[ \pi_1(x, y, z) = (x, y), \quad \pi_2(x, y) = x. \]
Furthermore, let $T_1$, $T_2$, and $T_3$ be the endomorphism on $X_1$, $X_2$, and $X_3$ represented by the matrices $\mathrm{diag}(2, 3, 4)$, $\mathrm{diag}(2, 3)$, and $\mathrm{diag}(2)$, respectively. Then $(X_i, T_i)_i$ and $(\pi_i)_i$ form a sequence of dynamical systems.

For a natural number $N$, denote by $Y_i|_N$ the restriction of $Y_i$ to its first $N$ coordinates, and let $p_{i, N}: Y_i|_N \rightarrow Y_{i+1}|_N$ be the projections for $i = 1, 2$. Since $Y_2$ and $Y_3$ are full shifts, we can use the same technique as in Example \ref{example: self affine sponges 1}. Therefore, we have for any exponent $\boldsymbol{a} = (a_1, a_2) \in [0, 1]^2$,
\[ h^{\boldsymbol{a}}(\boldsymbol{T}) = \lim_{N \to \infty} \frac{1}{N} \log{ \sum_{u \in \{0, 1\}^N} {\left( \sum_{v \in {p_{2, N}}^{-1}(u)} {| {p_{1, N}}^{-1}(v) |}^{a_1} \right)}^{a_2}}. \]

Now, let us evaluate $| {p_{1, N}}^{-1}(v) |$ using matrix products. This idea of using matrix products is due to Kenyon--Peres \cite{Kenyon--Peres: sofic}. Fix $(a, b) \in {\{0, 1\}}^2$ and let
\[ a_{ij} = | \{ e \in E \spa | \spa \text{$e$ is from $j$ to $i$ and the first two coordinates of its label is $(a, b)$} \} |. \]
Define a $3 \times 3$ matrix by $A_{(a, b)} = (a_{ij})_{ij}$. Then we have
\begin{equation*}
A_{(0, 0)} =
\begin{pmatrix}
0 & 1 & 1 \\
0 & 0 & 1 \\
1 & 1 & 0 \\
\end{pmatrix},
\spa A_{(0, 1)} =
\begin{pmatrix}
1 & 1 & 1 \\
1 & 1 & 0 \\
0 & 1 & 2 \\
\end{pmatrix},
\spa A_{(1, 0)} =
\begin{pmatrix}
1 & 2 & 2 \\
0 & 1 & 2 \\
2 & 2 & 1 \\
\end{pmatrix},
\spa A_{(1, 1)} = O.
\end{equation*}
Note that ${A_{(0, 0)}}^2 = A_{(0, 1)}$ and ${A_{(0, 0)}}^3 = A_{(1, 0)}$. For $v = (v_1, \cdots, v_N) \in Y_2|_N$ we have
\[ | {p_{1, N}}^{-1}(v) | \asymp \| A_{v_1} A_{v_2} \cdots A_{v_N} \|. \]
Here $A \asymp B $ means there is a constant $c > 0$ independent of $N$ with $c^{-1}B \leq A \leq cB$. For $\alpha = \frac{1 + \sqrt{5}}{2}$, we have $\alpha^2 = \alpha + 1$ and
\begin{equation*}
A_{(0, 0)}
\begin{pmatrix}
\alpha \\
1 \\
\alpha \\
\end{pmatrix}
= 
\begin{pmatrix}
1 + \alpha \\
\alpha \\
1 + \alpha \\
\end{pmatrix}
= \alpha
\begin{pmatrix}
\alpha \\
1 \\
\alpha \\
\end{pmatrix}, \quad
A_{(0, 1)}
\begin{pmatrix}
\alpha \\
1 \\
\alpha \\
\end{pmatrix}
= \alpha^2
\begin{pmatrix}
\alpha \\
1 \\
\alpha \\
\end{pmatrix}, \quad
A_{(1, 0)}
\begin{pmatrix}
\alpha \\
1 \\
\alpha \\
\end{pmatrix}
= \alpha^3
\begin{pmatrix}
\alpha \\
1 \\
\alpha \\
\end{pmatrix}.
\end{equation*}
Therefore,
\begin{equation*}
\| A_{v_1} A_{v_2} \cdots A_{v_N} \|
\asymp
\norm{
 A_{v_1} A_{v_2} \cdots A_{v_N}
\begin{pmatrix}
\alpha \\
1 \\
\alpha \\
\end{pmatrix}
} \asymp \lambda_{v_1} \lambda_{v_2} \cdots \lambda_{v_N}
\end{equation*}
where $\lambda_{(0,0)} = \alpha$, $\lambda_{(0,1)} = \alpha^2$, $\lambda_{(1,0)} = \alpha^3$.

Fix $u \in {\{0, 1\}}^{\mathbb{N}}$ and suppose there are $n$ numbers of zeros in $u$. Also, if there are $k$ numbers of $(0, 0)$s in $v = (v_1, \cdots, v_N) \in {p_{2, N}}^{-1}(u)$, there are $n - k$ numbers of $(0, 1)$s and $N - n$ numbers of $(1, 0)$s in $v$. Then
\[ {\lambda_{v_1}}^{a_1} \cdots {\lambda_{v_N}}^{a_1} = \alpha^{a_1k} \alpha^{2a_1(n-k)} \alpha^{3a_1(N-n)}. \]
Therefore,
\begin{align*}
\sum_{v \in {p_{2, N}}^{-1}(u)} {| {p_{1, N}}^{-1}(v) |}^{a_1}
&= \sum_{(v_1, \cdots, v_N) \in {p_{2, N}}^{-1}(u)} {\lambda_{v_1}}^{a_1} \cdots {\lambda_{v_N}}^{a_1} = \sum_{k=0}^n {n \choose k} \alpha^{a_1k} \alpha^{2a_1(n-k)} \alpha^{3a_1(N-n)} \\
&= {\left( \alpha^{a_1} + \alpha^{2a_1} \right)}^n \alpha^{3a_1(N-n)}.
\end{align*}
This implies
\begin{align*}
\sum_{u \in \{0, 1\}^N} {\left( \sum_{v \in {p_{2, N}}^{-1}(u)} {| {p_{1, N}}^{-1}(v) |}^{a_1} \right)}^{a_2}
&= \sum_{n=0}^N {N \choose n}{\left( \alpha^{a_1} + \alpha^{2a_1} \right)}^{a_2n} \alpha^{3a_1a_2(N-n)} \\
&= \left\{ {\left( \alpha^{a_1} + \alpha^{2a_1} \right)}^{a_2} + \alpha^{3a_1a_2} \right\}^N.
\end{align*}
We conclude that
\begin{align*}
h^{\boldsymbol{a}}(\boldsymbol{T})
&= \lim_{N \to \infty} \frac{1}{N} \log{\left\{ {\left( \alpha^{a_1} + \alpha^{2a_1} \right)}^{a_2} + \alpha^{3a_1a_2} \right\}^N} \\
&= \log{ \left\{ {\left( {\left( \frac{1 + \sqrt{5}}{2} \right)}^{a_1} + {\left( \frac{3 + \sqrt{5}}{2} \right)}^{a_1} \right)}^{a_2} + {(2+\sqrt{5})}^{a_1a_2} \right\}}.
\end{align*}

As in Example \ref{example: self-affine sponges two}, the Hausdorff dimension of $X_1$ is obtained by letting $a_1 = \log_4{3}$ and $a_2 = \log_3{2}$;
\begin{align*}
\mathrm{dim}_H(X_1)
&= \log{ \left\{ {\left( {\left( \frac{1 + \sqrt{5}}{2} \right)}^{\log_4{3}} + {\left( \frac{3 + \sqrt{5}}{2} \right)}^{\log_4{3}} \right)}^{\log_3{2}} + \sqrt{{(2+\sqrt{5})}} \right\}} \\
&= 1.4598\cdots.
\end{align*}

\section*{Acknowledgement}
I am deeply grateful to my mentor, Masaki Tsukamoto, who not only has reviewed this paper several times throughout the writing process but has patiently helped me understand ergodic theory in general with his expertise.

I also want to thank my family and friends for their unconditional support and everyone who has participated in my study for their time and willingness to share their knowledge. This work could not have been possible without their help.

\vspace{0.5cm}

\address{
Department of Mathematics, Kyoto University, Kyoto 606-8502, Japan}

\textit{E-mail}: \texttt{alibabaei.nima.28c@st.kyoto-u.ac.jp}

\end{document}